\documentclass[12pt]{article}
\usepackage{comment}
\input epsf
\usepackage{epsfig}
\usepackage{latexsym}
\usepackage{amsfonts}
\usepackage{amsmath}
\usepackage[psamsfonts]{amssymb}

\def\H{{\mathbf{\overline{H}}}}

\def\R{{\mathbf{R}}}
\def\S{{\mathbf{S}}}
\def\D{{\mathbf{D}}}
\def\P{{\mathbf{P}}}

\def\C{{\mathbf{C}}}

\def\P{\mathcal{P}}
\def\Q{\mathcal{Q}}

\def\Z{{\mathbf{Z}}}

\newtheorem{thm}{Theorem}[section]
\newtheorem{example}[thm]{Example}
\newtheorem{prop}[thm]{Proposition}
\newtheorem{lemma}[thm]{Lemma}

\newtheorem{cor}[thm]{Corollary}
\newtheorem{rmk}[thm]{Remark}
\newtheorem{df}[thm]{Definition} 

\numberwithin{equation}{section}

\author{Alexandre Eremenko\thanks{Supported by NSF grant DMS-1361836.},
Andrei Gabrielov\thanks{Supported by NSF grant DMS-1161629.}$\, $
and Vitaly Tarasov}
\title{Spherical quadrilaterals with three non-integer angles}
\begin{document}
\maketitle
\begin{abstract} A spherical quadrilateral is a bordered surface
homeomorphic to a closed disk, with four distinguished
boundary points called corners, equipped with a Riemannian
metric of constant curvature $1$, except at the corners, and such that
the boundary arcs between the corners are geodesic.
We discuss the problem of classification of these quadrilaterals
and perform the classification up to isometry in the case that
one corner of a quadrilateral is integer
(i.e., its angle is a multiple of $\pi$)
while the angles at its other three corners are not multiples of $\pi$.
The problem is equivalent to classification of Heun's equations with real
parameters and unitary monodromy, with the trivial monodromy at one of its
four singular point.


\end{abstract}

\noindent
\section{Introduction}\label{sec:Intro}

A {\em circular polygon}
is a simply connected surface whose boundary consists of
arcs of circles.

More precisely, let $\S$ be the Riemann sphere. A circle $C$
is the set of fixed points of an anti-conformal involution of $\S$.
This definition does not depend on the metric but only on the conformal
structure of $\S$. Let $\D$
be the closed unit disk, and $A=(a_0,a_1\ldots,a_{n-1})$ a sequence of
boundary points enumerated
according to the positive orientation of $\partial\D$. Let
$f:\D\to\S$ be a continuous function analytic in the open disk $U$, which is a
local homeomorphism on $\D\setminus A$, and such that
each {\em side} $(a_j,a_{j+1})\subset\partial\D$ is mapped by $f$
into some circle $C_j$; the circles are not necessarily distinct.
The triple
$(\D, A, f)$ is called a circular polygon, $f$ is called
the {\em developing map}, and the points of $A$ are called {\em corners}.

Two polygons $(\D,A,f)$ and $(\D,A',g)$ are {\em equivalent} if there
exist conformal maps $\phi:\S\to\S$ and $\psi:\D\to\D,\;
\psi(a_j)=a^\prime_j,\; 0\leq j\leq n-1,$ such that $\phi\circ f=g\circ\psi.$
Note that equivalence preserves the linear order of the corners, not only
their cyclic order.

At each corner $a_j$ the (interior) angle $\alpha_j\geq 0$ is defined.
We measure angles in multiples of $\pi$, so angle $\alpha_j$ contains $\pi\alpha_j$
radians.

If $\S$ is equipped with the standard spherical metric, and the circles
$C_j$ are geodesic (great circles), the polygon is called a
{\em spherical polygon}.

Two spherical polygons are considered equal (isometric) if they are equivalent,
and in addition, $\phi$ is an isometry of the sphere.
Notice that isometry of polygons with vertices $(a_0,\ldots,a_{n-1})$
and $(a_0^\prime,\ldots,a_{n-1}^\prime)$, according to our definition
must send $a_j$ to $a_j^\prime$.

One can give an alternative definition of a spherical polygon: it is a
surface homeomorphic to the closed disk, equipped with a Riemannian metric
of constant curvature $+1$ except at the corners $a_0,a_1,\ldots,a_{n-1}$
on the boundary, and such that the sides $(a_j,a_{j+1})$ are
geodesic, and the metric has conical singularities at the corners.
Two spherical polygons in the sense of this definition are considered
equal if there is an orientation-preserving isometry respecting
the labels of the corners.

The equivalence of the two definitions follows from the consideration
of the developing map: if $z_0$ is an interior point of the spherical polygon
in the sense of the second definition, there is a neighborhood of this point
and an isometry $f$ from this neighborhood to the sphere with the standard
spherical metric. Analytic continuation of $f$ gives the developing map.

Notice that great circles $C_j$ corresponding to a spherical polygon
satisfy the following
\vspace{.1in}

\noindent
{\bf Condition C.} {\em Every two circles either coincide or intersect
at two points transversally.}
\vspace{.1in}

This paper is a part of the series whose goal is classification of spherical
quadrilaterals up to isometry.
A complete classification of spherical $2$-gons is contained
in \cite{Troy1}, and of spherical triangles in \cite{E,FFKRUY}, see also
\cite{Klein-book}. Polygons with two non-integer angles were studied
in \cite{EGT1,EGT2}. In this paper we study spherical polygons, especially
quadrilaterals, with three non-integer angles. The general context and history
of the question are described in \cite{EGT2}. We also mention Schilling's
thesis \cite{Schil} where circular quadrilaterals with one angle $2\pi$
were studied in great detail.

It will be convenient to use the upper half-plane $\H$ instead of
the unit disk $\D$ for the parametrization of a polygon.
In that case, developing map $f$ of a circular $n$-gon satisfies a Schwarz
differential equation
\begin{equation}\label{schwarz1}
\frac{f'''}{f'}-\frac{3}{2}\left(\frac{f''}{f'}\right)=R,
\end{equation}
where
\begin{equation}\label{schwarz2}
R(z)=\sum_{j=0}^{n-1}\frac{1-\alpha_j^2}{2(z-a_j)^2}+\frac{c_j}{z-a_j},
\end{equation}
$a_j$ and $c_j$ are real, $\alpha_j\geq 0$, and $c_j$ satisfy three linear relations
$$\sum_{j=0}^{n-1}c_j=0,
\quad\sum_{j=0}^{n-1}a_jc_j=-\sum_{j=0}^{n-1}\frac{1-\alpha_j^2}{2},
\quad\sum_{j=0}^{n-1}a_j^2c_j=-\sum_{j=0}^{n-1}a_j(1-\alpha_j^2).$$
These three conditions ensure that the point at $\infty$ is not singular.
Equations (\ref{schwarz1}) and (\ref{schwarz2}) give
a parametrization of all circular polygons with
corners $a_j$ and angles $\alpha_j$. When $a_j$ and $\alpha_j$
are fixed, this set depends on $n-3$ real parameters,
which are called {\em accessory parameters}.
The condition that the polygon is equivalent to a spherical polygon translates
into the condition that the monodromy of the equation (\ref{schwarz1}) is
conjugate to a subgroup of the projective unitary group $PSU(2)$
(we say that the monodromy is {\em unitarizable}); this condition imposes $n-3$ equations
on the accessory parameters.

Thus our problem of classifying spherical polygons is equivalent to the study
of real solutions of these $n-3$ equations.
We mainly restrict ourselves to the case $n=4$ and three non-integer angles,
but in the next section we prove
a statement which holds for every $n$.

The position of four corners of a quadrilateral up to a conformal
automorphism of the disk depends on one parameter
which is called the {\em modulus}. One accessory parameter must be
determined from the condition that the monodromy is unitarizable.
So in general we can expect finitely many quadrilaterals with prescribed
angles and prescribed modulus. Our two main problems are, first, to
describe the necessary conditions on the angles for
a quadrilateral with the given angles to exist, for some value of the modulus;
and second, when those conditions are satisfied, to find upper and lower bounds
for the number of quadrilaterals with the fixed angles and modulus.

The plan of the paper is the following.

In section \ref{sec:condition} we discuss conditions that the angles
of a spherical quadrilateral with three non-integer angles must satisfy.
We also write explicitly the equation on the accessory parameter
for this case. This equation is algebraic and
we determine its degree.
This degree gives an upper estimate for the number
of quadrilaterals with prescribed modulus and angles.

In the following sections we classify these spherical quadrilaterals up to isometry,
and obtain a lower estimate for the number of quadrilaterals with prescribed
modulus and angles.

\section{Conditions on the angles and equation for the accessory parameter}\label{sec:condition}

\noindent
{\bf 1.} Consider a finite collection of circles on the sphere $\S$ satisfying Condition C,
such that their union is connected and there are no triple intersections.

These circles define a cell decomposition of the sphere whose
$0$-, $1$- and $2$-cells are called {\em vertices}, {\em edges} and {\em faces}.

To each pair $(F,v)$ where $F$ is a face and $v\in\partial F$ is a vertex
on the boundary of $F$ an angle $\alpha\in (0,1)$ is assigned.
It is the interior angle of the circular polygon $F$ at $v$,
measured in multiples of $\pi$.

Let us assign an arbitrary orientation to each circle.
Then to each pair $(F,v)$ we assign a sign $+$ or $-$
by the following rule: if the orientations of the two boundary edges of $F$
adjacent at $v$ are consistent we assign $+$, otherwise we assign $-$.

Let $f:\D\to \S$ be the developing map
of a circular polygon with corners $a_j$ whose sides are mapped to the circles
of our collection.

The $f$-preimage of the cell decomposition of the sphere is called a {\em net};
it is a cell decomposition of the closed disk.
The net is a combinatorial object: two nets are equivalent if there is
a homeomorphism of $\D$ which sends one to another and respects the
labels of the corners.


At every corner $v$ of the polygon we have an angle,
which is the sum of the angles
assigned to $(f(F),f(v))$ over all faces $F$ adjacent to $v$.
This angle is an integer when the two sides of the polygon adjacent to $v$ are mapped to the same
circle, and non-integer otherwise.

The signs assigned to faces $f(F)$ for $F$ adjacent to $v$ make an alternating
sequence. For a corner with non-integer angle the length of this sequence
is odd, and we assign to $v$ the sign $s(v)=+$ if the sequence
begins and ends with $+$, otherwise $s(v)=-$.

With this setting we have

\begin{prop}\label{prop1}
The number $N$ of signs $s(v)=-$ assigned to non-integer
corners is congruent modulo $2$ to
\begin{equation}\label{sigma0}
\sigma_0:=\sum_{ j:\alpha_j\not\in\Z}[\alpha_j]+
\sum_{ j:\alpha_j\in\Z}(\alpha_j-1).
\end{equation}
\end{prop}

{\em Proof.}
Consider a corner $v$ with integer angle $\alpha$.
Orientation of the sides changes at $v$ if and only if $\alpha-1$
is odd. So the number of changes of orientation at the corners with integer
angles is congruent modulo 2 to
$$\sum_{ j:\alpha_j\in\Z}(\alpha_j-1).$$
Now consider a corner $v$ with non-integer angle $\alpha$.
It is easy to see that the sides of the polygon change orientation
at $v$ if $[\alpha]$ is even and
$s(v)=-$, or if $[\alpha]$ is odd and $s(v)=+$.
In the other two cases the orientation of the sides does not change at $v$.

We conclude that the number
$$\sum_{ j:\alpha_j\in\Z}(\alpha_j-1)+\sum_{ j:\alpha_j\not\in\Z}[\alpha_j]+N$$
is congruent modulo $2$ to
the number of changes of orientation as we trace the boundary of the net,
but this last number is evidently even.
Thus the number of signs $s(v)=-$ assigned to the corners is congruent to
$\sigma_0$.
\vspace{.1in}

\begin{figure}
\centering
\includegraphics[width=4.0in]{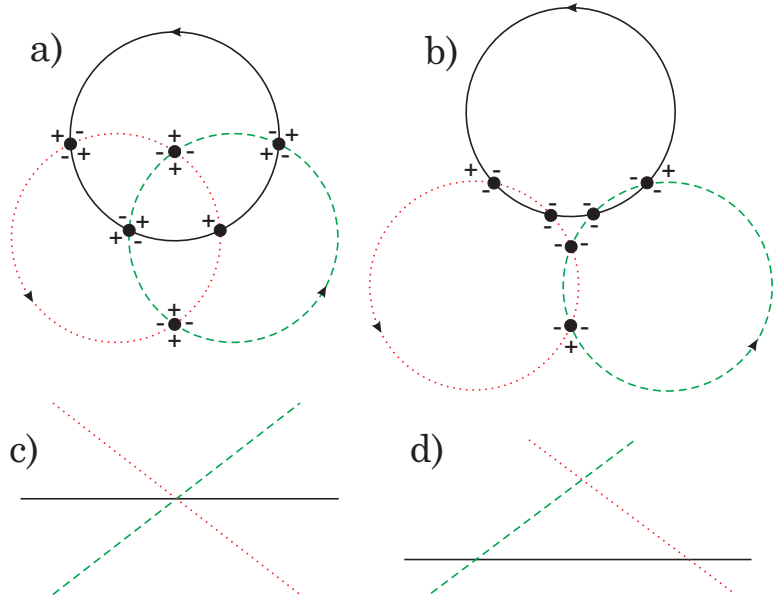}
\caption{ a) Generic unitarizable case, b) and d)
Non-unitarizable cases,\newline
c) Exceptional case.}\label{3circles-all}
\end{figure}

\noindent
{\bf 2.} Consider three circles on the sphere such
that every two of them intersect, and no two are tangent to each other.
There are four possible topologically distinct configurations shown
in Fig.~\ref{3circles-all}.

Indeed, if there is a triple intersection, we can send this point
to $\infty$ and obtain a configuration of three lines. There are two possible
configurations of three lines, shown in Fig.~\ref{3circles-all}c
and Fig.~\ref{3circles-all}d.

If there is no triple intersection, one of the circles $C_j$ either separates
the two intersection points of the other two circles or does not separate them.
If it does separate, then every circle separates the intersection points
of the other two circles. So we have two possibilities shown in
Fig.~\ref{3circles-all}a and Fig.~\ref{3circles-all}b.

It is easy to see that the triple of circles is equivalent to a geodesic
triple exactly in the cases a) and c). The case c) will be called
{\em exceptional}.

For a circular polygon whose sides are mapped by the developing map
into a generic configuration of three circles (that is the circles satisfy Condition C
and are not exceptional) we give a criterion in terms of the angles for the case a).

\begin{prop}\label{prop2}
Let $Q$ be a circular polygon with non-integer angles
$\theta$, $\theta'$, $\theta''$ and the rest of the angles integers.
Suppose that the images of its sides under the developing map
satisfy Condition C. Then $Q$ is equivalent to a spherical
polygon if and only if it is either exceptional, or
\begin{equation}\label{3}
\Sigma:=\cos^2\pi\theta+\cos^2\pi\theta'+\cos^2\pi\theta''+
2(-1)^\sigma
\cos\pi\theta\cos\pi\theta'\cos\pi\theta''<1,
\end{equation}
where $$\sigma=\sum_{ j:\alpha_j\in\Z}(\alpha_j-1).$$
\end{prop}

{\em Proof.}
Suppose that there are no triple intersections. Then Proposition \ref{prop1}
is applicable. Let us choose the orientation of the circles as shown
in Fig.~\ref{3circles-all}. This means that we choose a triangular
face in the cell decomposition,
orient its boundary counterclockwise, and this gives the orientation
of all three circles. Then we put signs according to the rule described above.
The original triangle gets signs $+,+,+$ on its corners.
Notice that each pair of circles intersects twice, and at each intersection
we have four angles, two of them with $+$ and two with $-$.
All angles formed by two circles which are marked $+$ are equal, and all
angles which are marked with $-$ are equal. The sum
of an angle marked with $+$ and an angle marked $-$ in the configuration
of two circles equals $1$.

Consider the functions
$$\phi^\pm(\beta,\beta',\beta'')
=\cos^2\pi \beta+\cos^2\pi\beta'+\cos^2\pi\beta''\pm 2\cos\pi\beta\cos\pi\beta'\cos\pi\beta''.$$
Notice that for any integers $m,n,k$ we have
$$\phi^\pm(\beta,\beta',\beta')=\phi^\pm(\beta+m,\beta'+n,\beta''+k)
\quad\mbox{if}\quad m+n+k\quad\mbox{is even},$$
and
$$\phi^\pm(\beta,\beta',\beta'')=\phi^\mp(\beta+m,\beta'+n,\beta''+k)\quad\mbox{if}\quad m+n+k\quad\mbox{is odd}.$$
Moreover,
$$\phi^\pm(\beta,\beta',\beta'')=\phi^\mp([\beta]+1-\{\beta\},\beta',\beta''),$$
and similarly for $\beta',\beta''$.

For the angles $\beta,\beta',\beta''$ of a triangular face in a)
we have $\phi^+(\beta,\beta',\beta'')<1$,
and for a triangular face
in b) we have $\phi^+(\beta,\beta',\beta'')>1$.
The angles of these triangles are in $(0,1)$. Each angle
$\theta,\theta',\theta''$ corresponds to a vertex of the cell
decomposition of the sphere where a pair of circles intersect.
These three pairs are all different. So, to obtain $\theta,\theta',\theta''$
from $\beta,\beta',\beta''$, we have to add the integer parts $[\theta]$,
$[\theta'], [\theta'']$ to $\beta,\beta',\beta''$, and replace some angles
by complementary angles. Notice that for corners marked with $-$ the fractional
part equals to an angle marked $-$ on the cell decomposition of the sphere.

Thus passing from $\beta,\beta',\beta''$ to $\theta,\theta',\theta''$
will require
$$[\theta]+[\theta']+[\theta'']+N$$
sign changes in $\phi^{\pm}$, where $N$ is the number of corners
marked with $-$. By Proposition \ref{prop1}, $N$ has the same parity
as $\sigma_0$. This proves Proposition \ref{prop2}.
\vspace{.2in}

{\bf 3.} The Schwarz equation (\ref{schwarz1}),
(\ref{schwarz2}) can be reduced to a linear differential equation
with regular singularities. In the case of quadrilaterals to which we restrict
from now on, this linear differential equation is the Heun equation (\ref{heun})
below.
Assuming that the quadrilateral is parametrized by
the upper half-plane, we can choose
the corners at the points
$0,\; 1,\; a$ and $\infty$. (The cyclic order of the corners depends on
the position of $a$ on the real line, so the notation $\alpha_0,\alpha_1$
differs in this section from the rest of the paper.) Then
$f=y_1/y_0$, the ratio of two linearly
independent solutions of the equation
\begin{equation}\label{heun}
z(z-1)(z-a)\left\{ y''+\left(\frac{1-\alpha_0}{z}+\frac{1-\alpha_1}{z-1}+
\frac{1-\alpha_a}{z-a}\right)y'\right\}+\alpha'\alpha''zy=\lambda y,
\end{equation}
where $\lambda$ is the accessory parameter,
\begin{equation}\label{alpha3}
\alpha'-\alpha''=\alpha_\infty,
\end{equation}
and
\begin{equation}\label{fuchs}
\alpha_0+\alpha_1+\alpha_a+\alpha'+\alpha''=2.
\end{equation}
We assume that $\alpha_0$ is an integer, while
$\alpha_1,\alpha_a,\alpha_\infty$
are not integers.

The exponents at $0$ are $0$ and $\alpha_0>0$. There is always one
holomorphic solution at $0$ which can be normalized
by $y_1(z)=z^{\alpha_0}+\ldots$. The second linearly independent solution
may in general contain logarithms. Condition C in the introduction is
equivalent to the absence of logarithms, and it must be satisfied if
the monodromy is unitarizable. To write this condition explicitly, we substitute
a power series $y_0(z)=1+c_1z+c_2z^2+\ldots$ to (\ref{heun})
and try to determine the coefficients. We obtain the recurrence relation
$$a(1-\alpha_0)c_1-\lambda=0,$$
\begin{equation}\label{rec}
r_nc_{n+1}-(q_n+\lambda)c_n+p_nc_{n-1}=0, \quad n\geq 1,
\end{equation}
where
$$r_n=(n+1)(n+1-\alpha_0)a$$
$$q_n=n\left((n-\alpha_0)(1+a)+a(1-\alpha_1)+1-\alpha_a\right),$$
$$p_n=(n-1+\alpha')(n-1+\alpha'').$$
As $\alpha_0$ is an integer, we have
$r_{\alpha_0-1}=0$, so $c_{\alpha_0}$ can be only determined when
$$G_{\alpha_0}(\lambda):= -(q_{\alpha_0-1}+\lambda)c_{\alpha_0-1}+
p_{\alpha_0-1}c_{\alpha_0-2}=0.$$
It is easy to see that $G_{\alpha_0}$ is a polynomial of degree $\alpha_0$.
Thus the equation
\begin{equation}\label{logs}
G_{\displaystyle\alpha_0}(\lambda)=0
\end{equation}
gives a necessary and sufficient condition for the absence of logarithms
in the solutions of the Heun equation (\ref{heun}).
Now we have

\begin{thm}\label{thm1}
Suppose that $\alpha_0\geq 2$ is a positive integer while the other
three
angles are not integers, and
$a\in\R$.

Then for the projective monodromy of (\ref{heun}) to be unitarizable, it
is necessary and sufficient that
(\ref{logs}) holds
and either we have (\ref{3}), or $\Sigma=1$ in (\ref{3}) and
we have the exceptional case (as in Fig.~\ref{3circles-all}c).
For every given angles, there are only
finitely many values of  parameter $a$ corresponding to the exceptional case.
\end{thm}

{\em Proof.} Suppose that the polygon defined by (\ref{heun})
is not exceptional. Necessity of conditions (\ref{logs}) and (\ref{3})
has already been established. For sufficiency, notice that all parameters
in (\ref{heun}) are real, so solutions with real initial conditions
at some non-singular point of the real line will map a real neighborhood
of this point into the real line. Thus any solution of the Schwarz equation
will map this neighborhood into a circle.
Condition (\ref{logs}) implies that two
adjacent sides cannot be tangent unless they coincide. As we assume that
the configuration is not exceptional, Proposition \ref{prop2} implies that the
configuration of the circles to which the sides belong is equivalent to
a configuration of great circles. Therefore the monodromy group,
which is generated by products of even numbers of reflections in these circles,
is conjugate to a subgroup of $PSU(2)$.

Let us prove now that the number of equivalence classes of
exceptional quadrilaterals with
prescribed angles is finite. Partitions of the sphere corresponding to such
quadrilaterals consist of $3$ intersecting lines and, possibly, one
vertex on one of these lines. These configurations with prescribed
angles between the lines are conformally equivalent. Let us show that there are
finitely many possible nets. The area of a spherical polygon is a function
of its interior angles, so when the angles are fixed, the area is fixed as well.
This gives an upper bound for the number of faces of the net,
so we have finitely many nets.
This completes the proof.
\vspace{.1in}

Spherical polygons correspond to {\em real} solutions of equation
(\ref{logs}). Thus for given angles and fixed modulus $a$, there exist
at most $\alpha_0$ spherical quadrilaterals with these angles and modulus.
In the rest of the paper we will prove a lower estimate.
This will be done by obtaining a combinatorial classification of all possible
quadrilaterals with given angles.

We also mention here an algebraic argument which gives a lower estimate.
Three-term recurrence (\ref{rec}) shows that (\ref{logs}) is the characteristic
equation of a finite Jacobi matrix
\begin{equation}\label{J}
J=\left(\begin{array}{cccccccc}
0&r_0&&&&&&\\
p_1&-q_1&r_1&&&&&\\
&p_2&-q_2&r_2&&&&\\
&&&&\ldots&&&\\
&&&&&p_{n-1}&-q_{n-1}&r_{n-1}\\
&&&&&&p_n&-q_n\end{array}\right).
\end{equation}
For the theory of such matrices we refer to \cite{Krein}.

If $r_jp_{j+1}>0$, $0\leq j\leq n-1$, then all eigenvalues are real and simple.
If $p_j\neq 0$ for $1\leq j\leq n$, then we have
\begin{equation}
\label{form}
J^TR=RJ,
\end{equation}
\def\diag{\mathrm{diag}\, }
where $R=\diag(s_0,\ldots,s_n)$, $s_0=1$ and
$$s_j=s_{j-1}\frac{r_{j-1}}{p_j}, \quad 1\leq j\leq n.$$
Suppose that the sequence $s_j$, $0\leq j\leq n$ has $P$ positive and $N$
negative terms. Then $J$ has at most $\min\{ P,N\}$
pairs of complex conjugate non-real
eigenvalues. The reason is that equation (\ref{form}) says that $J$ is
Hermitian with respect to the quadratic form $R$ of signature $(P,N)$,
and the estimate of the number of non-real eigenvalues follows
from Pontrjagin's theorem \cite{Pont}.

This gives an algorithm, which for any given angles and given $a$ provides
a lower estimate of the number of spherical
quadrilaterals with these parameters. This lower estimate is not always
exact.
Another lower estimate which will be proved in Section 7
is conjectured to be exact.
For example, the choice $\alpha_0=4,\alpha_1=\alpha_\infty=1/2,
\alpha_a=3/2$ and $a<0$ gives by this algorithm that equation
(\ref{logs}) has at most two non-real solutions, while the result
we obtain in Section 7 implies that all solutions are
real in this case.
\vspace{.1in}

We will need the following observation. If one sets $a=0$ in the
recurrence relations (\ref{rec}), we obtain
$r_j=0$, $0\leq j\leq n-1$, so the Jacobi matric $J$ is lower triangular.
The eigenvalues in this case are $-q_j=j(\alpha_0+\alpha_a)-j-1)$,
$0\leq j\leq n=\alpha_0-1$.
They are all real and distinct, because $\alpha_a$ is not integer.
It follows that for $a$ sufficiently small, all solutions of (\ref{logs})
are real. We state this as

\begin{prop}\label{allreal}
When the modulus of a quadrilateral is sufficiently
large or small, that is $a$ in (\ref{heun}) is sufficiently close to $0$ or $1$,
all solutions of (\ref{logs}) are real, and there exist exactly $\alpha_0$
quadrilaterals with given angles and given modulus.
\end{prop}

\section{Introduction to nets}\label{intro-nets}

A spherical $n$-gon $Q$ (or a {\em spherical polygon} when $n$ is not specified)
has its corners labeled $a_0,\dots,a_{n-1}$
in the counterclockwise order on the boundary of $Q$.
This defines linear order on the corners of $Q$.
When $n=2,\,3$ and $4$, we call $Q$ a spherical {\em digon}, {\em triangle}
and {\em quadrilateral}, respectively.
The angles at the corners of $Q$ are measured in multiples of $\pi$.
A corner is integer if it has an integer, i.e., multiple of $\pi$, angle.
For $n=1$, there is a unique 1-gon with the angle $\alpha_0=1$ at its single corner.
For convenience, we often drop ``spherical'' and refer simply to $n$-gons,
polygons, etc.

Let $Q$ be a spherical polygon and $f:Q\to\S$ its developing map.
The images of the sides $(a_j,a_{j+1})$ of $Q$ are contained in
geodesics (great circles) on $\S$. These circles define a {\em partition}
$\P$ of $\S$ into vertices (intersection points of the circles) edges
(arcs of the circles between the vertices) and faces (components of the complement
to the circles). Two sides of $Q$ meeting
at its corner are mapped by $f$ into the same circle
if and only if the corner is integer.

The {\em order} of a corner is the integer part of its angle.
A {\em removable} corner is an integer corner of order 1.
A polygon $Q$ with a removable corner is isometric to a polygon with a smaller number of corners.

A polygon with all integer corners is called {\em rational}.
(The developing map extends to a rational function in this case, which
explains the name).
All sides of a rational polygon map to the same circle.

\begin{df}\label{def:net}
{\rm Preimage of $\P$ defines a cell decomposition $\Q$ of $Q$, called the {\em net} of $Q$.
The corners of $Q$ are vertices of $\Q$.
In addition, $\Q$ may have {\em side vertices}
and {\em interior vertices}.
If the circles of $\P$ are in general position, interior vertices have degree 4,
and side vertices have degree 3.
Each face $F$ of $\Q$ maps one-to-one onto a face of $\P$.
An edge $e$ of $\Q$ maps either onto an edge of $\P$ or onto a part of an edge of $\P$.
The latter possibility may happen when $e$ has an end at an integer corner of $Q$.
The adjacency relations of the cells of $\Q$ are compatible with the adjacency
relations of their images in $\S$.
The net $\Q$ is completely defined by its 1-skeleton, a connected
planar graph. When it does not lead to confusion, we use the same
notation $\Q$ for that graph.

If $C$ is a circle of $\P$, then the intersection $\Q_C$ of $\Q$ with the preimage
of $C$ is called the $C$-{\em net} of $Q$.
Note that the intersection points of $\Q_C$ with preimages of other circles of $\P$
are vertices of $\Q_C$.
A $C$-{\em arc} of $Q$ (or simply an arc when $C$ is not specified) is
a non-trivial path $\gamma$ in the 1-skeleton of $\Q_C$
that may have a corner of $Q$ only as its endpoint.
If $\gamma$ is a subset of a side of $Q$ then it is a {\em boundary arc}.
Otherwise, it is an {\em interior arc}.
The {\em order} of an arc is the number of edges of $\Q$ in it.
An arc is {\em maximal} if it is not contained in a larger arc.
Each side $L$ of $\Q$ is a maximal boundary arc.
The {\em order} of $L$ is, accordingly, the number of edges of $\Q$ in $L$.}
\end{df}

\begin{df}\label{df:primitive}
{\rm We say that $Q$ is {\em reducible} if it
contains a proper sub-polygon $Q'$ such that the corners of $Q'$ are some (possibly, all)
of the corners of $Q$.
Otherwise, $Q$ is {\em irreducible}.
The net of a reducible polygon $Q$ contains an interior arc
with the ends at two distinct corners of $Q$.
We say that $Q$ is {\em primitive} if it is irreducible and its net does not contain
an interior arc that is a loop.}
\end{df}

\begin{df}\label{df:combequiv}
{\rm Two irreducible polygons $Q$ and $Q'$ are
{\em combinatorially equivalent} if there is a
homeomorphism $Q\to Q'$ preserving the labels of their corners
and mapping the corners of $Q$ to the corners of $Q'$,
and the net $\Q$ of $Q$ to the net $\Q'$ of $Q'$.

Two rational polygons $Q$ and $Q'$ with all sides mapped
to the same circle $C$ of $\P$ are
combinatorially equivalent if there is a
homeomorphism $Q\to Q'$ preserving the labels of their corners
and mapping the net $\Q_C$ of $Q$ to the net $\Q'_C$ of $Q'$.

If $Q$ and $Q'$ are reducible and represented
as the union of two polygons $Q_0$ and $Q_1$ (resp., $Q'_0$ and $Q'_1$)
glued together along their common side, then $Q$ and $Q'$ are
combinatorially equivalent when there is a
homeomorphism $Q\to Q'$ inducing combinatorial equivalence between $Q_0$ and $Q'_0$,
and between $Q_1$ and $Q'_1$.}
\end{df}

Thus an equivalence class of nets is a
combinatorial object. It is completely determined by the labeling of the corners
and the adjacency relations.
We'll call such an equivalence class ``a net'' when this would not lead to confusion.

Conversely, given labeling of the corners and a partition $\Q$
of a disk with the adjacency relations
compatible with the adjacency relations of $\P$,
a spherical polygon with the
net $\Q$ can be constructed by gluing together the cells of $\P$
according to the adjacency relations of $\Q$.

In what follows we classify all equivalence classes of nets in the case when $n=4$
(that is, when $Q$ is a quadrilateral) and $\P$ is defined by three transversal great circles
(that is, $Q$ has one integer corner and three non-integer corners).
In this case, the boundary of each 2-cell of the net $\Q$ of $Q$ consists
of three segments mapped to the arcs of distinct circles, with the vertices at the common
endpoints of each two segments and, possibly, at the integer corner of $Q$.

By the Uniformization Theorem, each spherical quadrilateral is conformally equivalent
to a closed disk with four marked points on the boundary.
Thus conformal class of a quadrilateral $Q$ depends on one parameter, the {\em modulus} of $Q$.
In section \ref{sec:chains} below
we will study whether for given permitted angles of a quadrilateral an arbitrary modulus
can be achieved. This will be done by the method of continuity,
and for this we'll need some facts about deformation of spherical quadrilaterals.

\medskip
\begin{figure}
\centering
\includegraphics[width=2.5in]{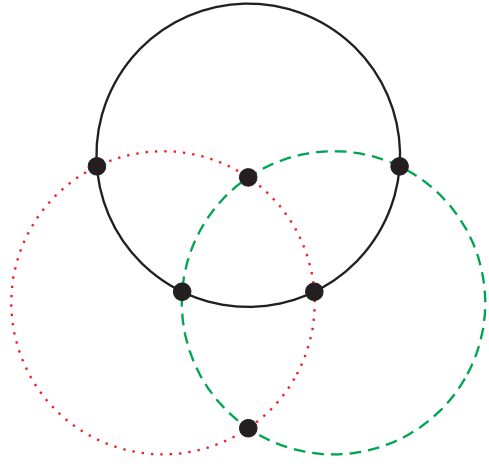}
\caption{Partition $\P$ of the Riemann sphere by three circles.}\label{partition}
\end{figure}

\section{Classification of primitive quadrilaterals}\label{sec:primitive}

Let us consider a partition $\P$ of the Riemann sphere $\S$
by three transversal great circles defined by the images of the sides
of a quadrilateral $Q$ with one integer corner (see Fig.~\ref{partition}).

We assume that with the corner $a_0$ of $Q$ selected so that the integer corner
of $Q$ is labeled $a_3$. Then the sides $L_3$ and $L_4$ of $Q$ are mapped
to the same circle $C$ of $\P$ (solid line in Fig.~\ref{partition}).
Let $C'$ and $C''$ be the circles of $\P$ to which the sides $L_1$ and $L_2$
of $Q$ are mapped (dotted and dashed lines in Fig.~\ref{partition}).

\begin{lemma}\label{corner}
The integer corner $a_3$ of an irreducible quadrilateral $Q$ does not map to a vertex of $\P$.
\end{lemma}

\noindent{\em Proof.}
Suppose that $a_3$ maps to the intersection of $C$ and $C'$ (the case of the intersection
of $C$ and $C''$ follows by symmetry).
Consider the union $G$ of all faces of $\Q$ adjacent to $a_3$.
Each of these faces is a triangle with two sides adjacent to $a_3$ mapped to $C$ and $C'$,
and a side opposite $a_3$ mapped to $C''$.
Thus $G$ is a triangle with an integer corner at $a_3$. Its base $\gamma$ is
an arc of $\Q$ mapped to $C''$.
Note that, since $Q$ is irreducible, all vertices of $\Q$ on $\gamma$ are interior vertices of $\Q$,
except the two endpoints of $\gamma$ which are on the boundary of $Q$.
It follows that, for each edge $e$ of $\gamma$, there is a face $F_e$ of $\Q$  adjacent to $e$ and not adjacent to $a_3$.
The union of all faces $F_e$ is a triangle $G'$ with the base $\gamma$ and a vertex $p$ common for all $F_e$.
Since $p$ is connected to $a_3$ by an interior arc, it cannot be a corner of $Q$.
This implies that the angle $\alpha_3$ of $Q$ at $a_3$ cannot be greater than $1$,
since $p$ cannot have degree greater than $4$.
If $\alpha_3=1$ then $p$ cannot be a boundary vertex of $\Q$, since then the whole boundary of $Q$ would map to $C$
which would contradict the assumption that $Q$ has three non-integer corners.
Also, $p$ cannot be an interior vertex of $\Q$, since then both ends of $\gamma$ would have degree at least $4$,
thus they should be both corners of $Q$, contradicting the assumption that $Q$ is irreducible.

\medskip
\begin{lemma}\label{no-interior-vertex}
The net $\Q$ of a primitive quadrilateral $Q$ does not have interior vertices
mapped to $C'\cap C''$.
\end{lemma}

\noindent{\em Proof.}
Let $q$ be an interior vertex of $\Q$ mapped to $C'\cap C''$.
Then there are four faces of $\Q$ adjacent to $q$.
The boundary $\Gamma$ of the union of those faces is a circle mapped to $C$.
Since $Q$ cannot be bounded by $\Gamma$, there is an interior edge of $\Q$ in $\Gamma$.
Let $\gamma$ be a maximal interior arc in $\Gamma$.
Since $Q$ is primitive, $\gamma$ has two distinct ends, which must be corners of $Q$
(the corner $a_1$ maps to $C'\cap C''$, and the corner $a_3$
cannot map to the intersection of two circles of $\P$ by Lemma \ref{corner}).
This contradicts the assumption that $Q$ is primitive.

\medskip
\begin{lemma}\label{interior-arcs}
The net $\Q$ of a primitive quadrilateral $Q$ does not have interior arcs of order greater
than 2. Any maximal interior arc of $\Q$ of order 2 which is mapped to $C$ has one of its 
ends at the corner $a_3$ of $Q$.
Any maximal interior arc of $\Q$ of order 1 which is mapped to $C$
has one of its ends at a corner of $Q$.
\end{lemma}

\noindent{\em Proof.}
We show first that a maximal interior arc $\gamma\subset\Q_C$ containing an interior vertex
$q$ of $\Q$ which is mapped to $C'$ has one end at the corner $a_3$ of $Q$.
If both edges of $\gamma$ adjacent to $q$ would have other ends
not at $a_3$ then, using the same argument as in the proof of Lemma \ref{no-interior-vertex},
one can show that the boundary $\Gamma$ of the union of the four faces of $\Q$ adjacent
to $q$ either would be a loop or would contain an interior arc connecting two corners of $Q$,
in contradiction with the assumption that $Q$ is primitive.
Thus $q$ is connected to $a_3$ by an edge of $\gamma$.
This implies that $\gamma$ has order 2, otherwise it would contain more than one
interior vertex of $\Q$, and only one such vertex can be connected to $a_3$.

Now we consider a maximal interior arc $\gamma$ of order 1 (i.e., a single edge of $\Q$) mapped to $C$.
If none of its ends would be at a corner of $Q$ then
each face $F$ of $\Q$ adjacent to $\gamma$ would have its sides other than $\gamma$
on the sides $L_1$ and $L_2$ of $Q$. Thus $F$ must have the corner $a_1$
as its vertex, but only one of the two faces of $\Q$ adjacent to $\gamma$ may be also adjacent to $a_1$.

\medskip
\begin{lemma}\label{corner-order-zero}
At least two non-integer corners of a primitive quadrilateral $Q$ have order zero.
\end{lemma}

\noindent{\em Proof.}
It is enough to show that if the corner $a_0$ of $Q$ has positive order then
its corners $a_1$ and $a_2$ have order zero. If $a_0$ has positive order then
there is a maximal interior arc $\gamma$ of $\Q_C$ and a maximal interior arc $\gamma'$ of $\Q_{C'}$,
each of them having one end at $a_0$, such that $\Q$ has a face $F$ adjacent to $L_1$ and $\gamma$,
and a face $F'$ adjacent to $\gamma$ and $\gamma'$.
Due to Lemma \ref{no-interior-vertex}, $\gamma'$ has order 1, and the other end $q'$ of $\gamma'$ is on $L_2$.
Thus $\gamma$ also have order 1, and its other end $q$ is on $L_2$.
This implies that $a_1$ has order zero, as it is a vertex of $F$.
This also implies that $a_2$ has order zero. Otherwise, the same arguments will show
that a maximal interior arc with one end at $a_2$ must have the other end on $L_1$,
which is impossible.

\begin{figure}
\centering
\includegraphics[width=5.0in]{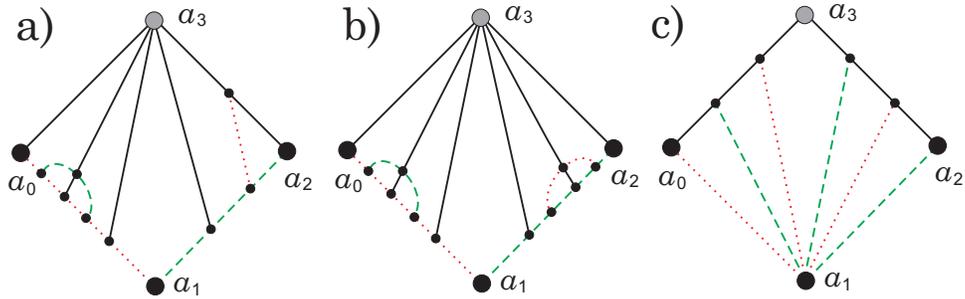}
\caption{Primitive quadrilaterals of types $X$ and $\bar X$.}\label{3circles-x}
\end{figure}

\begin{figure}
\centering
\includegraphics[width=3.4in]{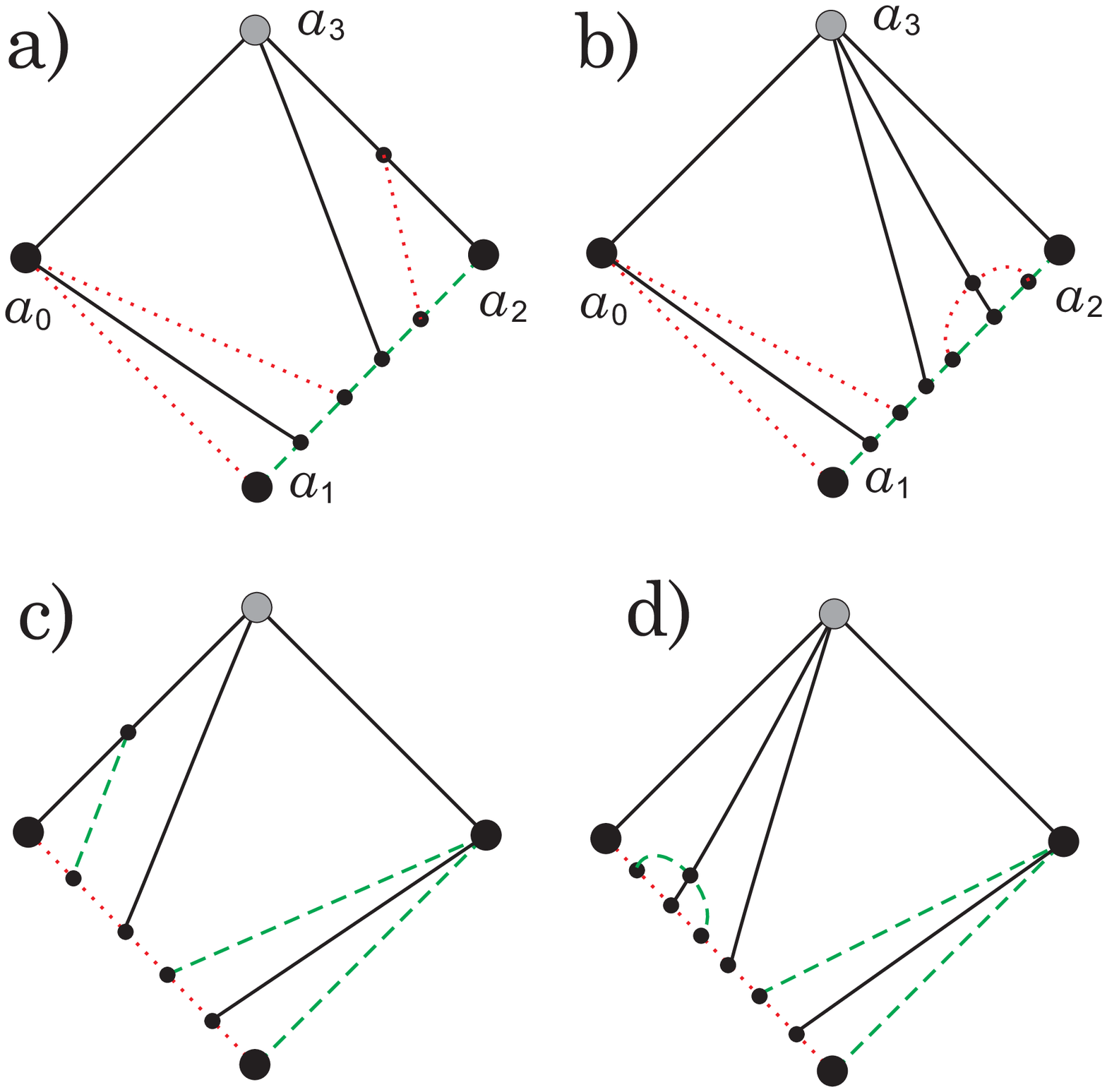}
\caption{Primitive quadrilaterals of types $R$ and $\bar R$.}\label{3circles-r}
\end{figure}

\begin{figure}
\centering
\includegraphics[width=3.4in]{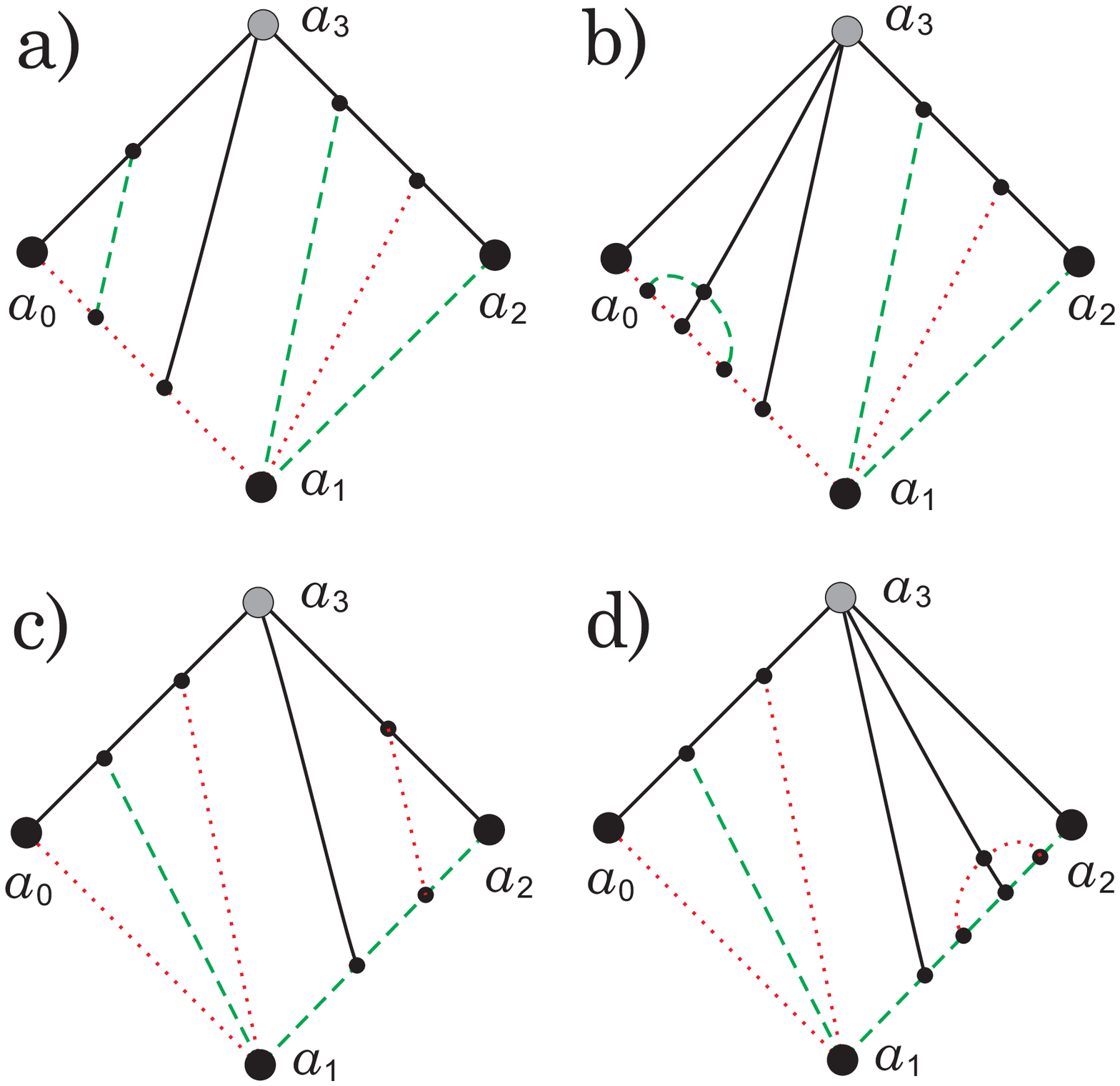}
\caption{Primitive quadrilaterals of types $U$ and $\bar U$.}\label{3circles-u}
\end{figure}

\begin{figure}
\centering
\includegraphics[width=5.0in]{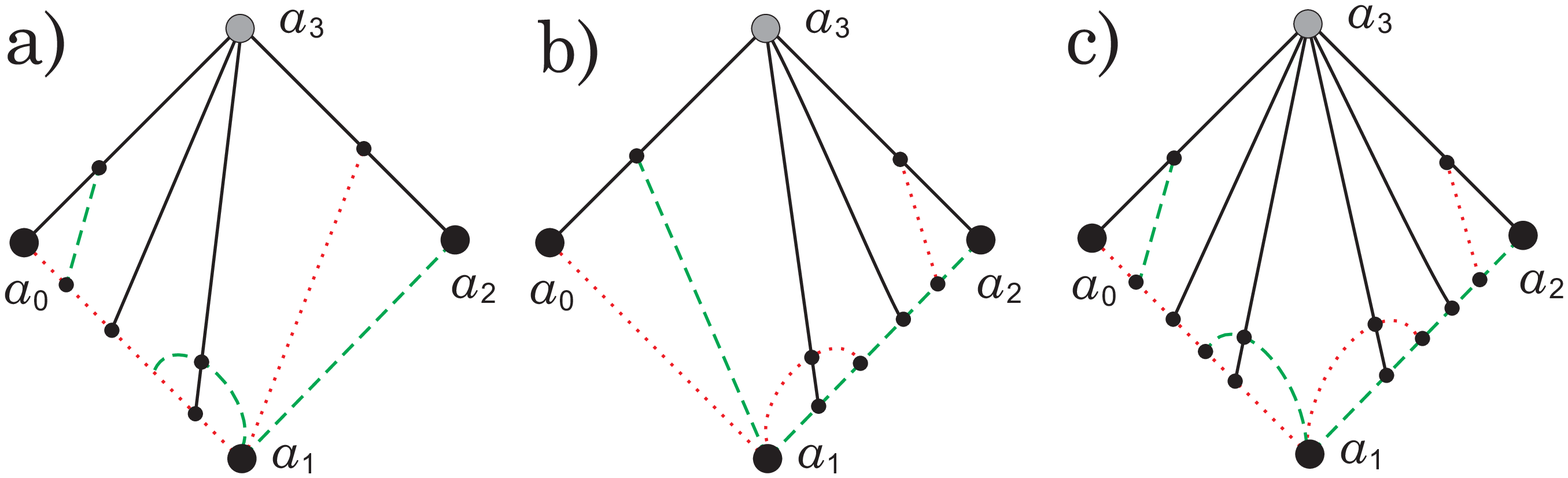}
\caption{Primitive quadrilaterals of types $V$, $\bar V$ and $Z$.}\label{3circles-vz}
\end{figure}

\medskip
\begin{thm}\label{thm:primitive}
Any primitive quadrilateral with one integer corner is equivalent to one of the following:\newline
$X_{\mu\nu}$ with $\mu\ge 0,\;\nu\ge 0$, the angles $[\alpha_0]=[\alpha_1]=[\alpha_2]=0,\;\alpha_3=\mu+\nu+1$,
$\mu$ (resp., $\nu$) interior arcs of the net having one end at $a_3$ and another end on $L_1$ (resp., $L_2$)
(Fig.~\ref{3circles-x}a and Fig.~\ref{3circles-x}b);\newline
$\bar X_{\mu\nu}$ with $\mu+\nu\ge 2$ even, the angles $[\alpha_0]=[\alpha_2]=[\alpha_3]=0,\;[\alpha_1]=(\mu+\nu)/2$,
$\mu$ (resp., $\nu$) interior arcs of the net having one end at $a_1$ and another end on $L_4$ (resp., $L_3$)
(Fig.~\ref{3circles-x}c);\newline
$R_{\mu\nu}$ with $\mu>0,\;\nu\ge 0$, the angles $[\alpha_0]=\mu,\;\alpha_3=\nu+1,\;[\alpha_1]=[\alpha_2]=0$,
$2\mu$ (resp., $\nu$) interior arcs of the net having one end at $a_0$ (resp., $a_3$) and another end on $L_2$
(Fig.~\ref{3circles-r}a and Fig.~\ref{3circles-r}b);\newline
$\bar R_{\mu\nu}$ with $\mu\ge 0,\;\nu>0$, the angles $\alpha_3=\mu+1,\;[\alpha_2]=\nu,\;[\alpha_0]=[\alpha_1]=0$,
$\mu$ (resp., $2\nu$) interior arcs of the net having one end at $a_3$ (resp., $a_2$) and another end on $L_1$
(Fig.~\ref{3circles-r}c and Fig.~\ref{3circles-r}d);\newline
$U_{\mu\nu}$ with $\mu>0,\;\nu>0$, the angles $\alpha_3=\mu+1,\;[\alpha_1]=\nu,\;[\alpha_0]=[\alpha_2]=0$,
$\mu$ (resp., $2\nu$) interior arcs of the net having one end at $a_3$ (resp., $a_1$) and another end on $L_1$ (resp., $L_3$)
(Fig.~\ref{3circles-u}a and Fig.~\ref{3circles-u}b);\newline
$\bar U_{\mu\nu}$ with $\mu>0,\;\nu>0$, the angles $[\alpha_1]=\mu,\;\alpha_3=\nu+1,\;[\alpha_0]=[\alpha_2]=0$,
$2\mu$ (resp., $\nu$) interior arcs of the net having one end at $a_3$ (resp., $a_1$) and another end on $L_4$ (resp., $L_2$)
(Fig.~\ref{3circles-u}c and Fig.\ref{3circles-u}d);\newline
$V_{\mu\nu}$ with $\mu>0,\;\nu>0$, the angles $\alpha_3=\mu+1,\;[\alpha_1]=\nu,\;[\alpha_0]=[\alpha_2]=0$, $\mu$
(resp., $2\nu-1$) interior arcs of the net having one end at $a_3$ (resp., $a_1$) and another end on $L_1$ (resp., $L_3$),
and one interior arc having one end at $a_1$ and another end on $L_1$ (Fig.~\ref{3circles-vz}a);\newline
$\bar V_{\mu\nu}$ with $\mu>0,\;\nu>0$, the angles $[\alpha_1]=\mu,\;\alpha_3=\nu+1,\;[\alpha_0]=[\alpha_2]=0$, $\nu$
(resp., $2\mu-1$) interior arcs of the net having one end at $a_3$ (resp., $a_1$) and another end on $L_2$ (resp., $L_4$),
and one interior arc having one end at $a_1$ and another end on $L_2$ (Fig.~\ref{3circles-vz}b);\newline
$\bar Z_{\mu\nu}$ with $\mu> 0,\;\nu> 0$, the angles $\alpha_3=\mu+\nu+1,\;[\alpha_0]=[\alpha_2]=0.\;[\alpha_1]=1$,
$\mu$ (resp., $\nu$) interior arcs of the net having one end at $a_3$ and another end on $L_1$ (resp., $L_2$),
one interior arc having one end at $a_1$ and another end on $L_1$, and
one interior arc having one end at $a_1$ and another end on $L_2$ (Fig.~\ref{3circles-vz}c).
\end{thm}

\begin{rmk}\label{arcs}
The nets of primitive quadrilaterals in Theorem \ref{thm:primitive}
may contain interior arcs with both ends on the sides of the quadrilateral
(see Figs.~\ref{3circles-x}-\ref{3circles-vz}). 
Location of these arcs is uniquely determined by the arcs of the net
having one end at a corner of the quadrilateral.
\end{rmk}

\noindent{\em Proof.}
If all non-integer corners have order $0$ then any maximal interior arc of $\Q$ mapped to $C$
has one end at $a_3$. Let $\mu$ (resp., $\nu$) be the number of such arcs having the other end on $L_1$
(resp, $L_2$). Then $Q$ is equivalent to $X_{\mu\nu}$.

If the corners $a_1$ and $a_2$ of $Q$ have order 0 but the order of $a_0$ is greater than 0,
then any maximal interior arc of $\Q$ mapped to $C$ has one end either at $a_0$ or at $a_3$.
Any such arc having one end at $a_0$ must have order 1, and its other end must be on $L_2$.
Also, any interior arc of $\Q$ mapped to $C'$ with one end at $a_0$ must have
order 1 and its other end on $L_2$, since it is the side of a face of $\Q$
adjacent to an interior arc mapped to $C$ with one end at $a_0$.
Since the order of $a_0$ is greater than 0, any maximal interior arc of $\Q$ mapped to $C$
with one end at $a_3$ must have its other end on $L_2$.
The interior arcs of $\Q$ (if any) having both ends on the sides of $Q$
but not at its corners are completely determined by the arcs having one end at a corner of $Q$.
If $\nu$ is the number of maximal interior arcs of $\Q$ with one end at $a_3$, and $\mu>0$ is the
number of interior arcs mapped to $C$ having one end at $a_0$ (which is equal to the number of interior
arcs mapped to $C'$ having one end at $a_0$) then $Q$ is equivalent to $R_{\mu\nu}$.
Similarly, if $a_0$ and $a_1$ have order 0 but the order of $a_2$ is greater than 0 then
$Q$ is equivalent to $\bar R_{\mu\nu}$ for some $\mu\ge 0$ and $\nu>0$.

If the corners $a_0$ and $a_2$ and $a_3$ of $Q$ have order 0 but the order of $a_1$ is greater than 0,
then any maximal interior arc of $Q$ having one end at $a_1$ has order 1 and its other end is either
on $L_3$ or $L_4$. The number of such arcs is even. There are no other arcs in $\Q$.
If $\mu$ (resp., $\nu$) is the number of interior arcs of $\Q$ having one end at $a_1$ and another end on $L_4$
(resp., $L_3$) then $Q$ is equivalent to $\bar X_{\mu\nu}$.

Consider now the case when $a_0$ and $a_2$ have order 0 but the orders of both $a_1$ and $a_3$ are greater
than 0. Suppose first that all maximal interior arcs of $\Q$ mapped to $C$ have one end at $a_3$
and the other end on $L_1$. Let $\mu>0$ be the number of these arcs.
Let $\nu>0$ be the number of maximal interior arcs of $\Q$ having one end at $a_1$ mapped to $C'$.
Then there are exactly $\nu$ maximal interior arcs with one end at $a_1$ mapped to $C''$.
There are two possibilities. If all these arcs have the other end on $L_3$ then $Q$ is equivalent to $U_{\mu\nu}$.
Alternatively, one of these arcs mapped to $C''$ may have the other end either on $L_1$ or (when $\mu=1$) on $L_4$.
Then $Q$ is equivalent to $V_{\mu\nu}$.
The interior arcs of $\Q$ (if any) having both ends on the sides of $Q$
but not at its corners are completely determined by the arcs having one end at a corner of $Q$.
Similarly, if all maximal interior arcs of $\Q$ mapped to $C$ have one end at $a_3$ and another end on $L_2$
then $Q$ is equivalent to either $\bar U_{\mu\nu}$ or $\bar V_{\mu\nu}$ for some $\mu>0$ and $\nu>0$.

Finally, suppose that $\mu>0$ (resp., $\nu>0$) of maximal interior arcs mapped to $C$ have one end at $a_3$
and the other end on $L_1$ (resp., $L_2$). Then there are exactly two interior arcs having one end at $a_1$,
one of them, mapped to $C'$, has the other end either on $L_2$ or (if $\nu=1$) on $L_3$, and the other arc,
mapped to $C''$, has the other end either on $L_1$ or (if $\mu=1$) on $L_4$. In that case $Q$ is equivalent
to $Z_{\mu\nu}$.

\begin{rmk}\label{primitive-equivalent} Sometimes it will be convenient to use notation $\bar X_{0,0}$, $R_{0,0}$, $\bar R_{0,0}$,
$U_{0,0}$, $\bar U_{0,0}$ for the quadrilaterals equivalent to $X_{0,0}$.
We may also use notation $U_{\mu\nu}$ and $\bar U_{\mu\nu}$ with either $\mu=0$ or $\nu=0$.
Then $U_{\mu,0}$ is equivalent to $X_{\mu,0}$, and $U_{0,\nu}$ is equivalent to $\bar X_{0,2\nu}$ for $\nu>0$.
Similarly, $\bar U_{0,\nu}$ is equivalent to $X_{0,\nu}$, and $\bar U_{\mu,0}$ is equivalent to $\bar X_{2\mu,0}$ for $\mu>0$.
Also, $V_{0,\nu}$ is equivalent to $\bar X_{1,2\nu-1}$, $\bar V_{\mu,0}$ is equivalent to $\bar X_{2\mu-1,1}$.
\end{rmk}

\section{Classification of irreducible quadrilaterals}\label{sec:irreducible}

\begin{lemma}\label{loop}
Let $\gamma$ be a loop of the net $\Q$ of an irreducible quadrilateral $Q$
with three non-integer corners.
Then $\gamma$ contains either the integer corner $a_3$
or the non-integer corner $a_1$ of $Q$.
\end{lemma}

\begin{figure}
\centering
\includegraphics[width=4.0in]{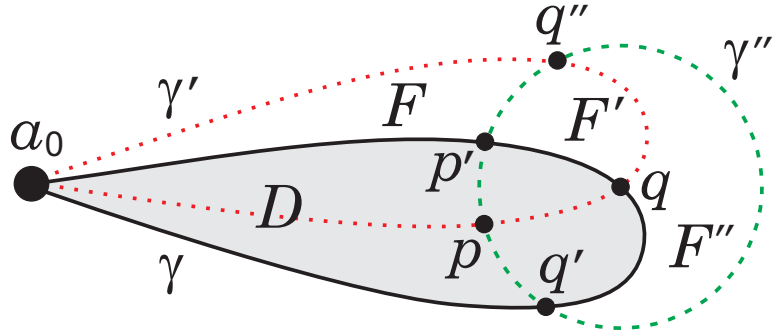}
\caption{Illustration of the proof of Lemma \ref{loop}.}\label{3circles-loop}
\end{figure}

\noindent{\em Proof.} A loop $\gamma$ bounds a disk $D\subset Q$
(shaded area of Fig.~\ref{3circles-loop})
consisting of 4 faces of $\Q$ with a common vertex $p$.
If $\gamma$ would not contain any corner then the 4 faces of $\Q$
adjacent to $\gamma$ outside $D$ would also have a common vertex,
thus $Q$ would be a sphere, a contradiction.
Thus $\gamma$ contains a corner of $Q$.
If $\gamma$ contains $a_0$ then it maps either to $C$ or to $C'$.
We assume that $\gamma$ maps to $C$, the other case being treated similarly.
Then $D$ contains part of a maximal arc $\gamma'$ of $\Q$ with one end at $a_1$
that maps to $C'$, and part of a maximal arc $\gamma''$ that maps to $C''$
(see Fig.~\ref{3circles-loop}).

We want to show that $\gamma'$ is a loop. Since $\gamma'$ contains $a_0$,
and $Q$ is irreducible, $\gamma'$ cannot contain any other corner of $Q$.
Let $q\ne a_0$ be the intersection point of $\gamma'$ and $\gamma$ (see Fig.~\ref{3circles-loop}).
Then the union of faces $F$ and $F'$ of $\Q$ adjacent to
the segment $[a_0,q]$ of $\gamma$ outside $D$ have a segment of $\gamma'$ outside $D$
connecting $q$ with $a_1$ as part of its boundary.

Next, we want to show that $\gamma''$ is a loop.
Let $p'$ be the intersection point of $\gamma''$ with $\gamma$ inside the disk bounded by $\gamma'$,
and let $q'\ne p'$ and $q''\ne p$ be the intersection points of $\gamma''$ with $\gamma$ and $\gamma'$,
respectively (see Fig.~\ref{3circles-loop}).
Then the face $F''$ of $\Q$ adjacent to the segments $[q,q']$ of $\gamma$ and $[q,q'']$ of $\gamma'$
has a segment $[q',q'']$ of $\gamma''$ in its boundary.
As $[q',q'']$ is a segment of the boundary of a single face $F''$ mapped to $C''$,
it cannot contain a corner of $Q$. Thus it is a single edge of $\Q$.
The boundary of $F''$ cannot contain a corner of $Q$,
since all vertices of $\Q$ on this boundary belong either to $\gamma$
or to $\gamma'$, and those two loops already contain the corner $a_0$.
Thus the loop $\gamma''$ contains no corners, a contradiction.

It follows that $\gamma$ cannot contain $a_0$.
A similar argument shows that $\gamma$ cannot contain $a_2$.
Thus $\gamma$ contains either $a_1$ or $a_3$.

\begin{figure}
\centering
\includegraphics[width=4.0in]{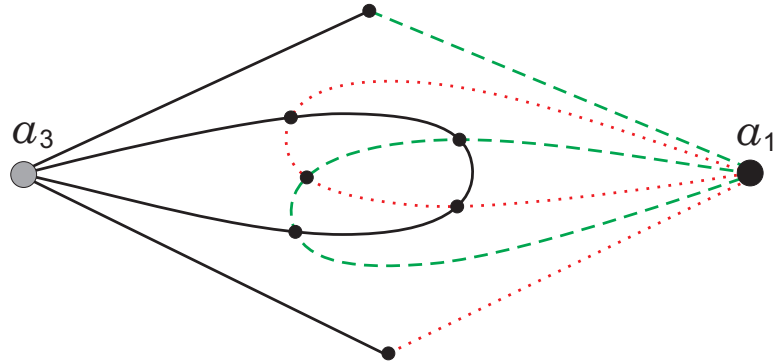}
\caption{Pseudo-diagonal connecting corners $a_1$ and $a_3$ of $Q$.}\label{3circles-pseudo-diagonal}
\end{figure}

\begin{prop}\label{pseudo}
Every loop of the net $\Q$ of an irreducible quadrilateral $Q$
with three non-integer corners is part of a pseudo-diagonal
(Fig.~\ref{3circles-pseudo-diagonal}) connecting
its integer corner $a_3$ with its non-integer corner $a_1$.
\end{prop}

\noindent{\em Proof.}
The same arguments as in the proof of Lemma \ref{loop} show that a loop containing $a_1$
should be part of a pair $(\gamma',\gamma'')$ where $\gamma'$ is a loop mapping
to $C'$ and $\gamma''$ is a loop mapping to $C''$, both loops containing $a_1$.
In addition, there should be a loop $\gamma$ mapped to $C$
intersecting both loops $\gamma'$ and $\gamma''$.
The loop $\gamma$ must contain a corner of $Q$, and the integer corner $a_3$
is the only possibility (any other corner is rejected by the same argument as in
the proof of Lemma \ref{loop}).

Conversely, if $\gamma$ is a loop of $\Q$ mapped to $C$ and containing $a_3$
then it bounds a disk $D\subset Q$ so that the arcs intersecting inside $D$
are parts of two loops, $\gamma'$ mapped to $C'$ and $\gamma''$ mapped to $C''$,
both containing $a_1$. This can be shown by the arguments similar to those in the
proof of Lemma \ref{loop}.

\medskip
\begin{cor}\label{irreducible}
Each irreducible spherical quadrilateral $Q$ with three non-integer corners
is either primitive or
can be obtained from a primitive polygon $\tilde Q$ of any type except $R_{\mu\nu}$ and $\bar R_{\mu\nu}$
with a face $F$ of its net containing the corners $a_1$ and $a_3$ of $\tilde Q$ are on its
 boundary, by adding any number of pseudo-diagonals connecting  $a_1$ with $a_3$ inside $F$.
\end{cor}

\medskip
We use notation $U_{\mu\nu}^\kappa$, $\bar U_{\mu\nu}^\kappa$,
$X_{\mu\nu}^\kappa$, $\bar X_{\mu\nu}^\kappa$, $V_{\mu\nu}^\kappa$, $\bar V_{\mu\nu}^\kappa$,
$Z_{\mu\nu}^\kappa$ for an irreducible quadrilateral obtained from
one of the primitive quadrilaterals in Theorem \ref{thm:primitive}
by adding $\kappa$ pseudo-diagonals.

\section{Classification of reducible quadrilaterals}\label{sec:reducible}

\begin{figure}
\centering
\includegraphics[width=5.0in]{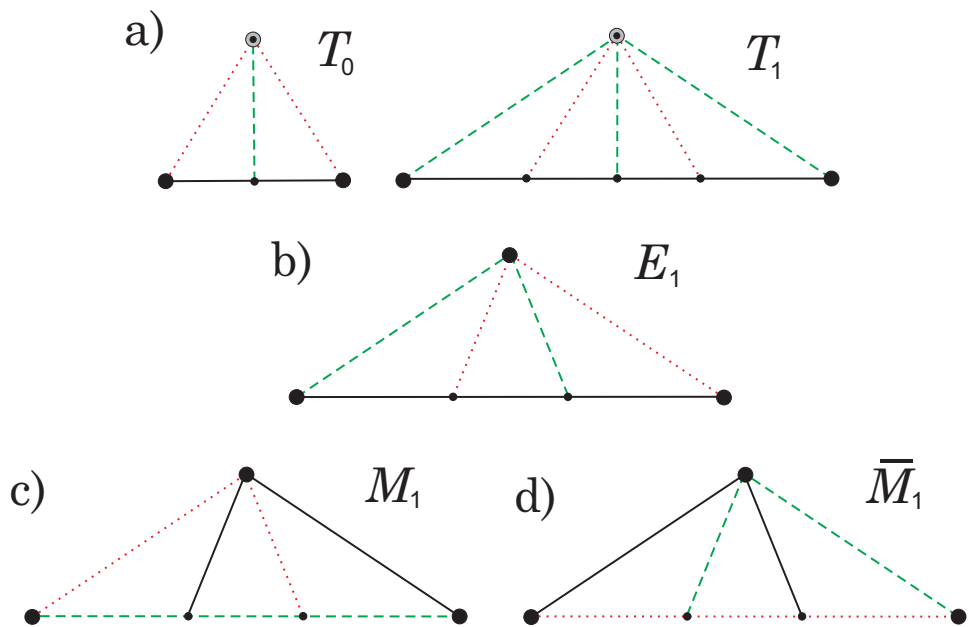}
\caption{Primitive triangles that appear in reducible quadrilaterals.}\label{3circles-tem}
\end{figure}

\begin{figure}
\centering
\includegraphics[width=5.0in]{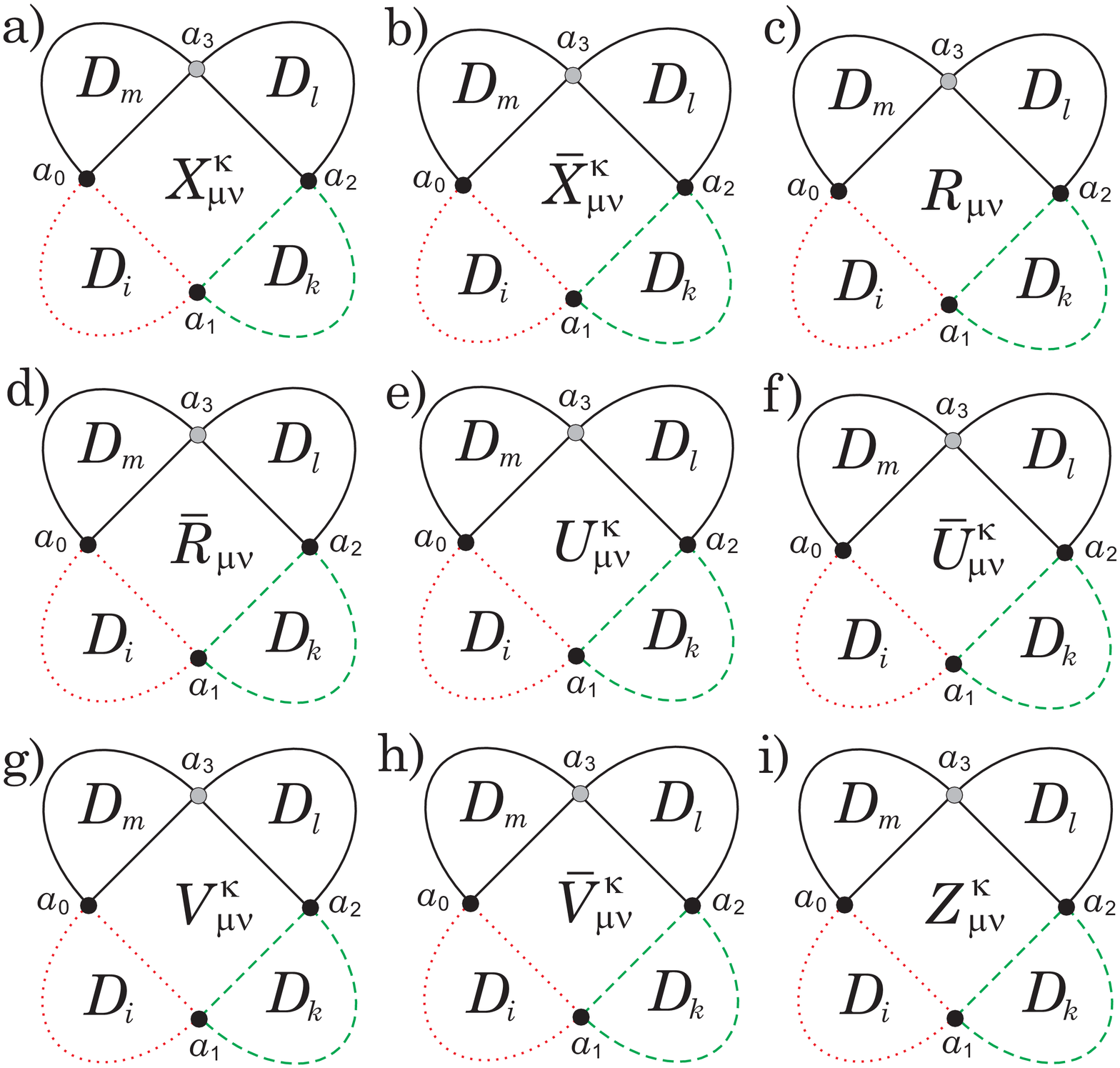}
\caption{Reducible quadrilaterals with a generic integer corner.}\label{3circles-xruvz}
\end{figure}

\begin{figure}
\centering
\includegraphics[width=5.0in]{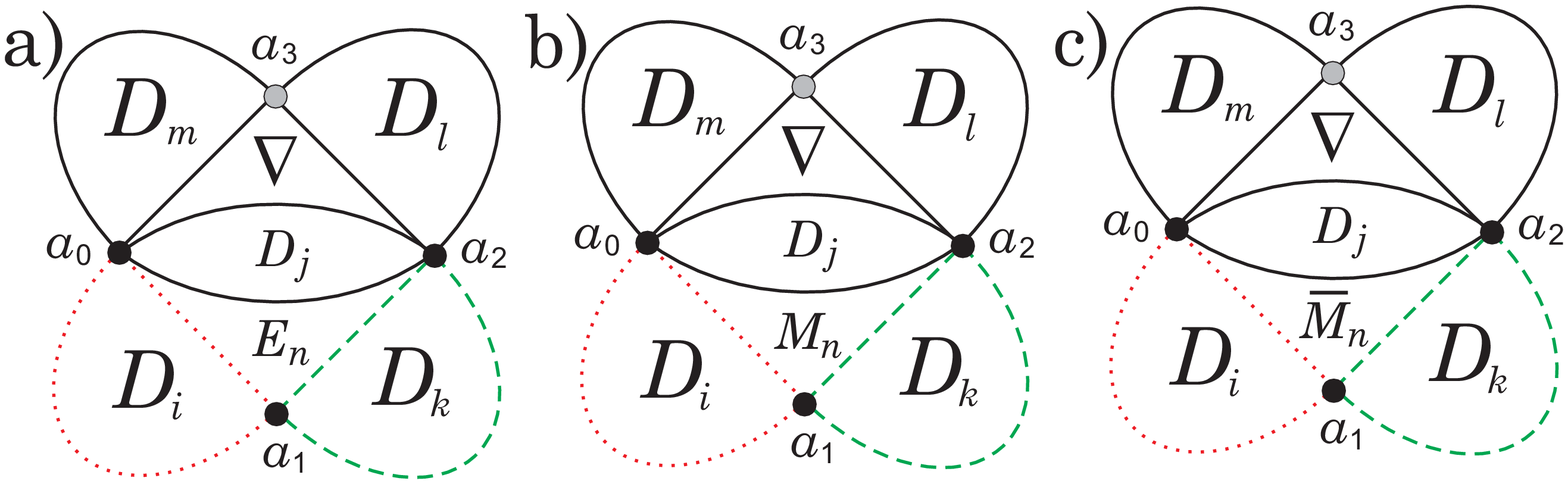}
\caption{Reducible quadrilaterals with a generic integer corner (cont).}\label{3circles-em}
\end{figure}

\begin{figure}
\centering
\includegraphics[width=5.0in]{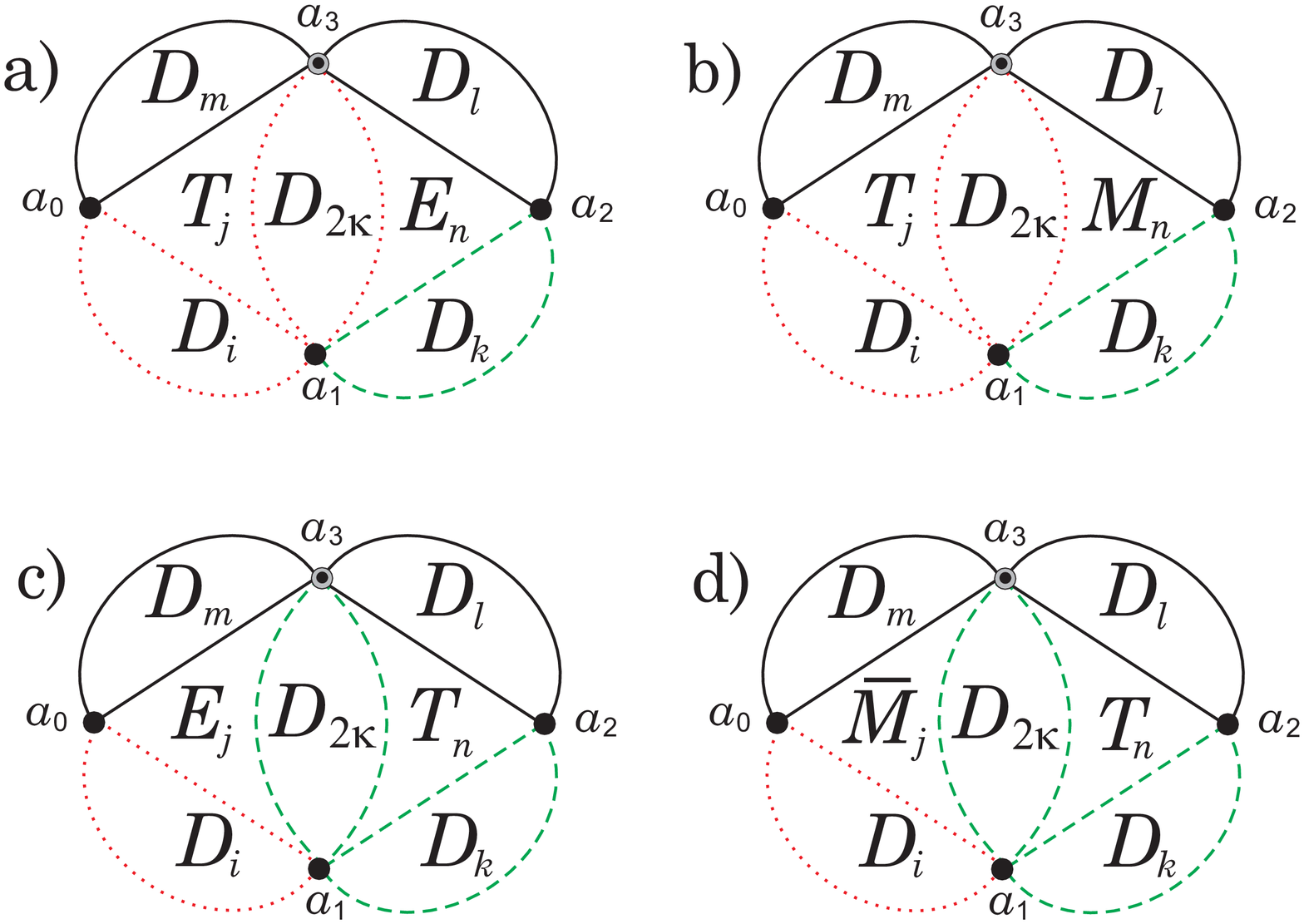}
\caption{Reducible quadrilaterals with a non-generic integer corner.}\label{3circles-tetm}
\end{figure}

\begin{figure}
\centering
\includegraphics[width=5.0in]{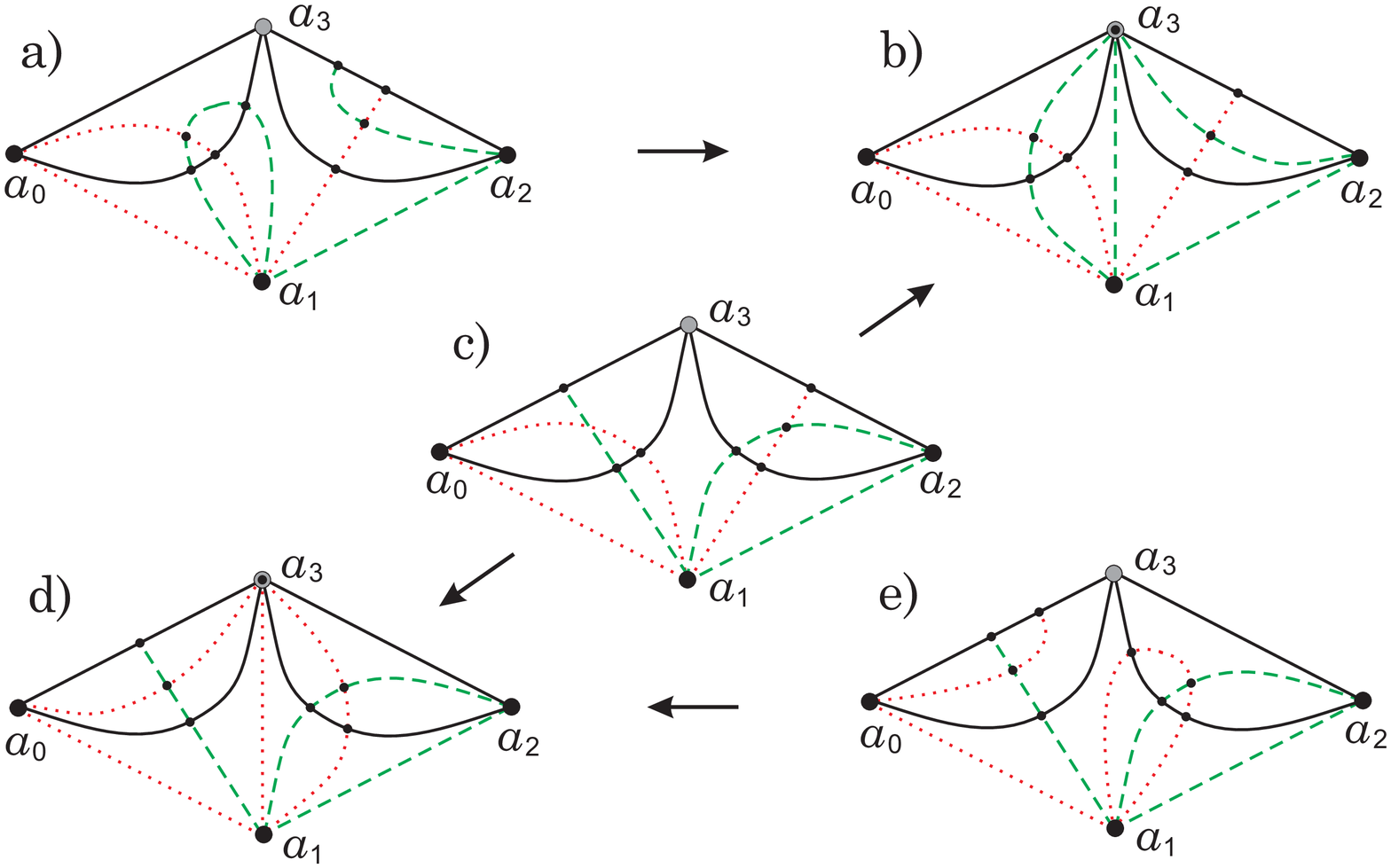}
\caption{An ab-chain of quadrilaterals of length 2.}\label{3circle-chain1}
\end{figure}

\begin{figure}
\centering
\includegraphics[width=5.0in]{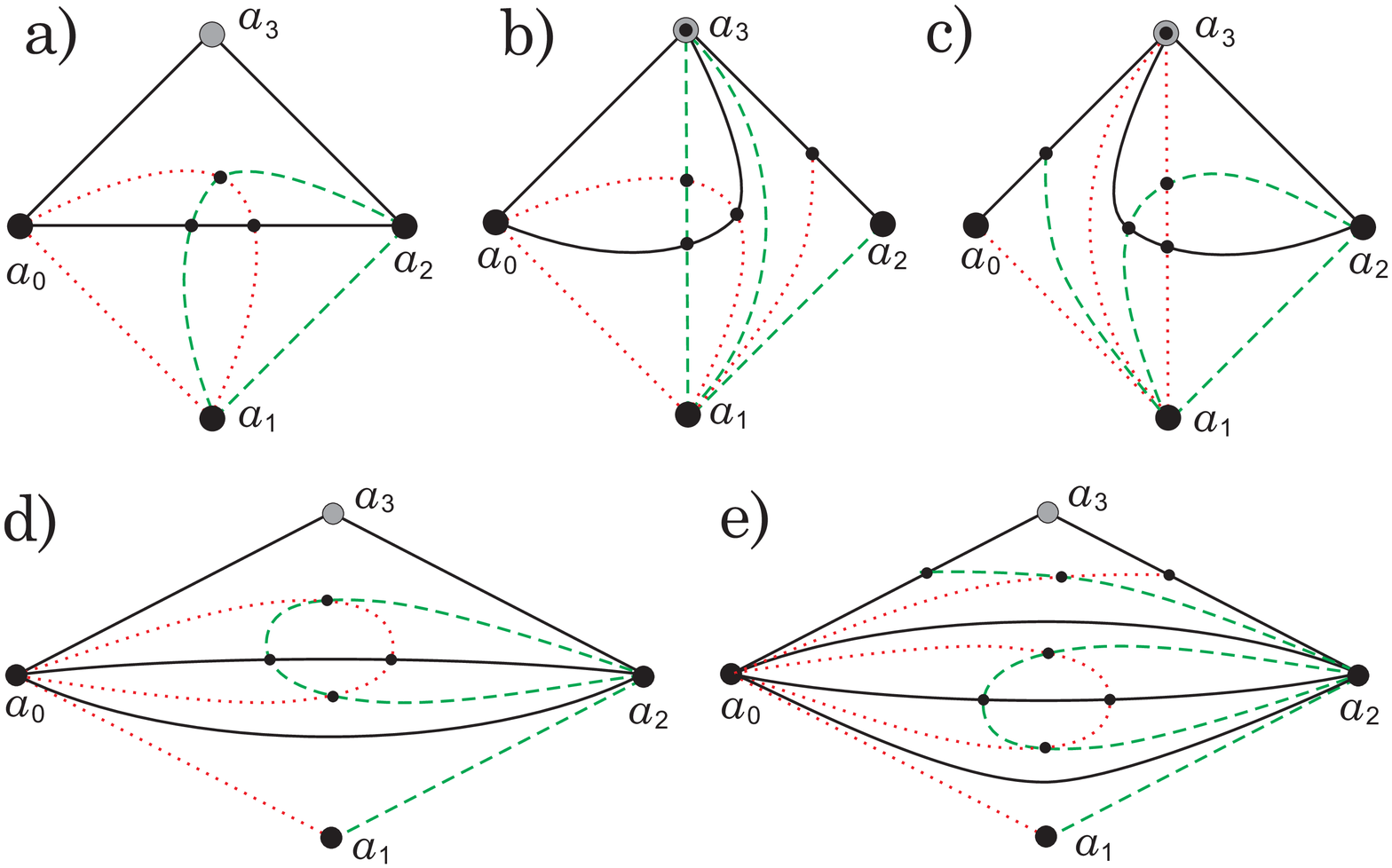}
\caption{Examples of reducible quadrilaterals.}\label{3circle-nabla}
\end{figure}

\medskip
To classify reducible quadrilaterals, we need digons $D_n$, and primitive triangles $T_n$,
$E_n$, $M_n$ and $\bar M_n$ (see Fig.~\ref{3circles-tem}).
We also need a triangle $\nabla$ with
all its angles equal 1 and all its sides mapped to the circle $C$ corresponding to the sides $L_3$ and $L_4$ of $Q$.
Geometrically, $\nabla$ is just a disk bounded by $C$ with three distinct points on $C$ marked as its corners.

A digon $D_n$ is a union of $n$ disks, each bounded by the same circle of $\P$,
glued together over the segments of their boundaries, each of them having its ends at the same two
points of the circle, so that their union is homeomorphic to a disk.
The angles at the two corners of $D_n$ are equal $n$.
The two boundary arcs of $D_n$ are called its sides.
We use the same notation $D_n$ for a digon having its sides on any of the three circles of $\P$,
and for any possible sizes of its two sides.
We use notation $D_0$ for an empty digon.

\medskip
A primitive triangle $T_n$ (see Fig.~\ref{3circles-tem}a)
has an integer corner with the angle $n\ge 1$
mapped to an intersection point of the two circles $C'$ and $C''$ of $\P$ corresponding
to the sides $L_1$ and $L_2$ of $Q$.
The side of $T_n$ opposite its integer corner maps to the circle $C$
of $\P$ corresponding to the sides $L_3$ and $L_4$ of $Q$.
We use the same notation $T_n$ for a triangle with two sides
mapped either to $C'$ or to $C''$.

A primitive triangle $E_n$ (see Fig.~\ref{3circles-tem}b)
has a non-integer corner with the integer part of the angle equal $n$
mapped to an intersection point of the two circles $C'$ and $C''$ of $\P$ corresponding
to the sides $L_1$ and $L_2$ of $Q$.
The side of $E_n$ opposite its integer corner maps to the circle $C$
of $\P$ corresponding to the sides $L_3$ and $L_4$ of $Q$.
The cyclic order of the sides of $E_n$ is consistent with the cyclic order of the
sides of $Q$.

A primitive triangle $M_n$ (resp., $\bar M_n$, see Fig.~\ref{3circles-tem}c
and Fig.~\ref{3circles-tem}d)
has a non-integer corner with the integer part of the angle equal $n\ge 1$
mapped to an intersection point of the circle $C'$ of $\P$ corresponding
to the side $L_1$ of $Q$ (resp., the circle $C''$ corresponding to the side $L_2$) with the circle $C$
corresponding to the sides $L_3$ and $L_4$.
The cyclic order of the sides of $M_n$ and $\bar M_n$ is consistent with the cyclic order of the
sides of $Q$.

If there is no arc of the net $\Q$ of a reducible quadrilateral $Q$
connecting its opposite corners then $Q$ contains an irreducible quadrilateral $\tilde Q$, with digons
$D_i$, $D_k$, $D_l$, $D_m$ attached to its four sides (see Fig.~\ref{3circles-xruvz}).
A non-empty disk may be attached only to a ``short'' side of $\tilde Q$ (that is, to a side
shorter than the full circle):\newline
If $\tilde Q=X_{\mu\nu}^\kappa$ then $i>0$ only if $\mu\le 1$ and $k>0$ only if $\nu\le 1$.\newline
If $\tilde Q=\bar X_{\mu\nu}^\kappa$ then $m>0$ only if $\mu\le 3$ and $l>0$ only if $\nu\le 3$.\newline
If $\tilde Q=R_{\mu\nu}$ then $k>0$ only if $\mu=1$ and $\nu=0$.\newline
If $\tilde Q=\bar R_{\mu\nu}$ then $i>0$ only if $\mu=0$ and $\nu=1$.\newline
If $\tilde Q=U_{\mu\nu}^\kappa$ then $i>0$ only if $\mu=1$ and $l>0$ only if $\nu=1$.\newline
If $\tilde Q=\bar U_{\mu\nu}^\kappa$ then $m>0$ only if $\mu=1$ and $k>0$ only if $\nu=1$.\newline
If $\tilde Q=V_{\mu\nu}^\kappa$ then $i>0$ only if $\mu=1$ and $l>0$ only if $\nu\le 2$.\newline
If $\tilde Q=\bar V_{\mu\nu}^\kappa$ then $m>0$ only if $\mu\le 2$ and $k>0$ only if $\nu=1$.\newline
If $\tilde Q=Z_{\mu\nu}^\kappa$ then $i>0$ only if $\mu=1$ and $k>0$ only if $\nu=1$.

If the corners $a_0$ and $a_2$ of $Q$ are connected by an arc of $\Q$
then $Q$ is a union of a triangle $\nabla$ and one of the
triangles $E_n$ (for $0\le 1\le 1$), $M_n$ (for $n\ge 1$), or $\bar M_n$ (for $n\ge 1$)
with a digon $D_j$ inserted between them and digons $D_i$, $D_k$, $D_l$, $D_m$ attached
to their other sides (see Fig.~\ref{3circles-em}).
Examples of such quadrilaterals are given in Fig.~\ref{3circle-nabla}a,d,e.

Note that the union of $\nabla$ and $E_1$ (see Fig.~\ref{3circle-nabla}a)
contains also arcs of $\Q$ connecting $a_1$ with $a_0$ and with $a_2$, and
is equivalent to one of the quadrilaterals either in Fig.~\ref{3circles-xruvz}c,d
or in Fig.~\ref{3circles-em}b,c.

The integer corner $a_3$ of $Q$ is ``generic'' -- it may map to a point on $C$ that
either belongs or does not belong to another circle of $\P$.
A non-empty disk may be attached only to a ``short'' side of a triangle
$M_n$ or $\bar M_n$:\newline
If $n>1$ then $k=0$ in Fig.~\ref{3circles-em}b.\newline
If $n>1$ then $i=0$ in Fig.~\ref{3circles-em}c.

If the corners $a_1$ and $a_3$ of $Q$ are connected by an arc of $\Q$
then $Q$ is a union of two primitive triangles (see Fig.~\ref{3circles-tetm})
with a digon $D_{2\kappa}$ inserted between them and and digons $D_i$, $D_k$, $D_l$, $D_m$
attached to their other sides.

We assume $n\ge 1$ in Fig.~\ref{3circles-tetm}b, and $j\ge 1$ in Fig.~\ref{3circles-tetm}d.
The integer corner $a_3$ of such a quadrilateral $Q$ is ``non-generic'' --
it must map to an intersection of $C$ with another circle of $\P$.
A non-empty disk may be attached only to a ``short'' side of any triangle:\newline
In Fig.~\ref{3circles-tetm}a, $m=0$ if $j>1$, and $l=0$ if $n>1$.\newline
In Fig.~\ref{3circles-tetm}b, $m=0$ if $j>1$, and $k=0$ if $n>1$.\newline
In Fig.~\ref{3circles-tetm}c, $m=0$ if $j>1$, and $l=0$ if $n>1$.\newline
In Fig.~\ref{3circles-tetm}d, $i=0$ if $j>1$, and $l=0$ if $n>1$.

Note that a digon with the odd angle cannot be inserted
between the two triangles in Fig.~\ref{3circles-tetm} since both sides of the digon
should be arcs of $\Q$ of length 1. One could consider also the union of triangles $M_1$ and $T_n$,
(or $T_j$ and $\bar M_1$) with a digon $D_{2\kappa+1}$ inserted between them (see Fig.~\ref{3circle-nabla}b and Fig.~\ref{3circle-nabla}c),
but the resulting quadrilateral would be equivalent to one of the quadrilaterals in Fig.~\ref{3circles-tetm}
with a digon $D_{2\kappa}$ inserted between two triangles.

\begin{thm}\label{integer-parts}
For any non-negative integers $A_0,\,A_1,\,A_2$ and a positive integer $A_3$ there exists a spherical quadrilateral
with $[a_0]=A_0$, $[a_1]=A_1$, $[a_2]=A_2$, and $a_3=A_3$.
\end{thm}

\noindent{\em Proof.}
Consider first the case $A_3=1$.
If $A_1\ge A_0+A_2$ then a quadrilateral $\bar X^0_{\mu\nu}$ with $\mu+\nu=2(A_1-A_0-A_2)$ and digons
$D_i$ and $D_k$ attached to its sides $L_1$ and $L_2$ (see Fig.~\ref{3circles-xruvz}b) where $i=A_0$
and $k=A_2$ has the required integer parts of its angles.
If $A_0+A_1+A_2$ is even and $(A_0-1,A_1,A_2-1)$ satisfy the triangle inequality then the union of
$\nabla$ and $E_0$ with a digon $D_j$ inserted between them and digons $D_i$ and $D_k$
attached to its sides $L_1$ and $L_2$ (see Fig.~\ref{3circles-tetm}a) where $j=(A_0+A_2-A_1)/2-1$,
$i=(A_0+A_1-A_2)/2$, $k=(A_1+A_2-A_0)/2$ has the required integer parts of its angles.
If $A_0+A_1+A_2$ is odd and $(A_0-1,A_1-1,A_2-1)$ satisfy the triangle inequality then the union of
$\nabla$ and $E_1$ with a digon $D_j$ inserted between them and digons $D_i$ and $D_k$
attached to its sides $L_1$ and $L_2$ (see Fig.~\ref{3circles-tetm}a) where $j=(A_0+A_2-A_1-1)/2$,
$i=(A_0+A_1-A_2-1)/2$, $k=(A_1+A_2-A_0-1)/2$ has the required integer parts of its angles.
If either $A_0\ge A_1+A_2+1$ or $A_2\ge A_0+A_1+1$ then
the union of $\nabla$ and either $M_n$ or $\bar M_n$ with $D_j$ inserted between them
and either $D_i$ attached to $L_1$ or $D_k$ attached $L_2$ (see Fig.~\ref{3circles-tetm}a and Fig.~\ref{3circles-tetm}b)
has the required integer parts of its angles.

Let now $A_3\ge 2$. If $A_1+A_3\le A_0+A_2$ then, attaching digons $D_l$ and $D_m$
to the sides $L_3$ and $L_4$ of the union of $\nabla$ and either $E_n$ or $M_n$ or $\bar M_n$,
for some $l$ and $m$ such that $l+m=A_3-1$,
we can reduce the problem to the case when $A_3=1$ and $A_1< A_0+A_2$.

Suppose now that $A_1+A_3=A_0+A_2+1$. Attaching digons $D_i$, $D_k$, $D_l$, $D_m$
to the sides of the quadrilateral $X_{0,0}$ (see Fig.~\ref{3circles-xruvz}a)
we can get a quadrilateral with the required integer parts of the angles.
Moreover, reducing (resp., increasing) by 1 the indices $i$ and $l$, and increasing
(resp., reducing) by 1 the indices $k$ and $m$, we obtain a quadrilateral with the
same angles. Thus we may always assume that either $\min(i,l)=0$ or $\min(k,m)=0$.
If, for example, $i=0$ or $k=0$, we may replace $X_{0,0}$ by either $X_{\mu,0}$
or $X_{0,\nu}$ and increase $A_3$. Similarly, if $m=0$ or $l=0$, we may replace
$X_{0,0}$ by either $\bar X_{\mu,0}$ or $\bar X_{0,\nu}$ and increase $A_1$.
It is easy to check that this way we can obtain all remaining values of the
integer parts of the angles.

\medskip
\begin{rmk} The primitive quadrilateral $\tilde Q$ contained in a reducible quadrilateral $Q$ is not unique.
The quadrilateral $Q$ shown in Fig.~\ref{3circle-chain1}c can be represented as one of the following:\newline
$X^0_{11}$ with $i=k=1$ and $l=m=0$ (see Fig.~\ref{3circles-xruvz}a);\newline
$\bar X^0_{22}$ with $i=k=0$ and $l=m=1$ (see Fig.~\ref{3circles-xruvz}b);\newline
$U^0_{11}$ with $i=l=1$ and $k=m=0$ (see Fig.~\ref{3circles-xruvz}e);\newline
$\bar U^0_{11}$ with $i=k=0$ and $l=m=1$ (see Fig.~\ref{3circles-xruvz}f).\newline
The quadrilateral $Q$ in Fig.~\ref{3circle-nabla}a is the union of a triangle $\nabla$
and a triangle $E_1$. The same quadrilateral can be represented as a quadrilateral
$R_{10}$ with $i=l=m=0$ and $k=1$, or as a quadrilateral $\bar R_{01}$ with $i=1$
and $k=l=m=0$ (see Fig.~\ref{3circles-xruvz}c,d).
\end{rmk}

\begin{lemma}\label{equivalent-reducible}
The following reducible quadrilaterals are equivalent:\newline
1) $X^\kappa_{\mu,1}$ with digons $D_i$, $D_k$, $D_l$, $D_m$ attached as in Fig.~\ref{3circles-xruvz}a, $k>0$,
and $U^\kappa_{\mu,1}$ with digons $D_i$, $D_{k-1}$, $D_{l+1}$, $D_m$ attached as in Fig.~\ref{3circles-xruvz}e.\newline
2) $X^\kappa_{1,\nu}$ with digons $D_i$, $D_k$, $D_l$, $D_m$ attached as in Fig.~\ref{3circles-xruvz}a, $i>0$,
and $\bar U^\kappa_{1,\nu}$ with digons $D_{i-1}$, $D_k$, $D_l$, $D_{m+1}$ attached as in Fig.~\ref{3circles-xruvz}f.\newline
3) $\bar X^\kappa_{2\mu,2}$ with digons $D_i$, $D_k$, $D_l$, $D_m$ attached as in Fig.~\ref{3circles-xruvz}b, $l>0$,
and $\bar U^\kappa_{\mu,1}$ with digons $D_i$, $D_{k+1}$, $D_{l-1}$, $D_m$ attached as in Fig.~\ref{3circles-xruvz}f.\newline
4) $\bar X^\kappa_{2,2\nu}$ with digons $D_i$, $D_k$, $D_l$, $D_m$ attached as in Fig.~\ref{3circles-xruvz}b, $m>0$,
and $U^\kappa_{1,\nu}$ with digons $D_{i+1}$, $D_k$, $D_l$, $D_{m-1}$ attached as in Fig.~\ref{3circles-xruvz}e.\newline
5) $\bar X^\kappa_{2\mu-1,3}$ with digons $D_i$, $D_k$, $D_l$, $D_m$ attached as in Fig.~\ref{3circles-xruvz}b, $l>0$,
and $\bar V^\kappa_{\mu,1}$ with digons $D_i$, $D_{k+1}$, $D_{l-1}$, $D_m$ attached as in Fig.~\ref{3circles-xruvz}h.\newline
6) $\bar X^\kappa_{3,2\nu-1}$ with digons $D_i$, $D_k$, $D_l$, $D_m$ attached as in Fig.~\ref{3circles-xruvz}b, $m>0$,
and $V^\kappa_{1,\nu}$ with digons $D_{i+1}$, $D_k$, $D_l$, $D_{m-1}$ attached as in Fig.~\ref{3circles-xruvz}g.\newline
7) $Z^\kappa_{\mu,1}$ with digons $D_i$, $D_k$, $D_l$, $D_m$ attached as in Fig.~\ref{3circles-xruvz}i, $k>0$,
and $V^\kappa_{\mu,1}$ with digons $D_i$, $D_{k+1}$, $D_{l-1}$, $D_m$ attached as in Fig.~\ref{3circles-xruvz}g.\newline
8) $Z^\kappa_{1,\nu}$ with digons $D_i$, $D_k$, $D_l$, $D_m$ attached as in Fig.~\ref{3circles-xruvz}i, $i>0$,
and $\bar V^\kappa_{1,\nu}$ with digons $D_{i-1}$, $D_k$, $D_l$, $D_{m+1}$ attached as in Fig.~\ref{3circles-xruvz}h.\newline
9) The union of $\nabla$ and $E_1$ with a digon $D_j$, $j>0$, inserted between them and digons $D_i$, $D_k$, $D_l$, $D_m$
attached as in Fig.~\ref{3circles-em}a, and either the union of $\nabla$ and $M_1$ with a digon $D_{j-1}$ inserted between them
and digons $D_i$, $D_{k+1}$, $D_l$, $D_m$ attached as in
Fig.~\ref{3circles-em}b,
or the union of $\nabla$ and $\bar M_1$ with a digon $D_{j-1}$ inserted between them
and digons $D_{i+1}$, $D_k$, $D_l$, $D_m$ attached as in
Fig.~\ref{3circles-em}c.\newline
10) The union of $\nabla$ and $E_1$ digons $D_i$, $D_k$, $D_l$, $D_m$ attached as in Fig.~\ref{3circles-em}a, $j=0$,
and either $R_{1,0}$ with digons $D_i$, $D_{k+1}$, $D_l$, $D_m$ attached as in Fig.~\ref{3circles-xruvz}c,
or $\bar R_{0,1}$ with digons $D_{i+1}$, $D_k$, $D_l$, $D_m$ attached as in Fig.~\ref{3circles-xruvz}d.
\end{lemma}

This can be shown by comparing the nets of the quadrilaterals in each of the ten cases.

\section{Chains of quadrilaterals}\label{sec:chains}

We use the results on deformation of quadrilaterals from \cite{EGT2}.

Let $\Gamma$ be
a family of curves in some region $D\subset\C$. Let $\lambda\geq 0$ be a measurable
function in $D$. We define the $\lambda$-length of a curve $\gamma$ by
$$\ell_\lambda(\gamma)=\int_\gamma\lambda(z)|dz|,$$
if the integral exists, and $\ell_\lambda(\gamma)=+\infty$ otherwise.
Then we set
$$L_\lambda(\Gamma)=\inf_{\gamma\in\Gamma}\ell_\lambda(\gamma),$$
and
$$A_\lambda(D)=\int_D\lambda^2(z)dm,$$
where $dm$ is the Euclidean area element.
Then the extremal length of $\Gamma$ is defined as
$$L(\Gamma)=\inf_\lambda\frac{L_\lambda^2(\Gamma)}{A_\lambda(D)}.$$
The extremal length is a conformal invariant. Extremal distance between two
closed sets is defined as the extremal length of the family of all curves in $D$
that connect these two sets.
For a flat rectangle with vertices $a_0$, $a_1$, $a_2$, $a_3$
conformally equivalent to the quadrilateral $Q$, the extremal distance between its
sides $[a_1,a_2]$ and $[a_3,a_0]$ is equal to $|[a_0,a_1]|/|[a_1,a_2]|$ \cite{Ahlfors}.

In addition to the extremal distance, we consider the ordinary intrinsic distances
between the pairs of opposite sides. They are defined as the infima of spherical lengths
of curves contained in our quadrilateral and connecting the two sides of a pair.

\begin{lemma}\label{modulus} (\cite[Lemma 15.1]{EGT2}
Consider a sequence of spherical quadrilaterals
whose developing maps $f$ are at most $p$-valent with a fixed integer $p$.
If the intrinsic distance between the sides $[a_1,a_2]$ and $[a_3,a_0]$ is bounded from below,
while the intrinsic distance between the sides $[a_0,a_1]$ and $[a_2,a_3]$ tends to zero,
then the extremal distance between the sides $[a_1,a_2]$ and $[a_3,a_0]$ tends to $+\infty$.
\end{lemma}

Let $Q$ be a quadrilateral with non-integer corners $a_0$, $a_1$ and $a_2$,
and with integer corner $a_3$ which is mapped to a point $X\in C$ that is not a vertex of $\P$.
Let $\Q$ be the net of $Q$.
When $X$ approaches a vertex $X_0$ of $\P$, and the combinatorial class
of $Q$ is fixed, we have the following possibilities.
\begin{itemize}
\item[(a)] $Q$ degenerates so that the intrinsic distance between its opposite sides
$L_1$ and $L_3$ tends to zero, while the intrinsic distance between its sides $L_2$ and $L_4$ does not tend to zero.
This happens when $\Q$ has an arc of order $1$ connecting $a_3$ with a
point on $L_1$ mapped to $X_0$,
but does not have an arc of order $1$ connecting $a_3$ to a point on $L_2$ mapped to $X_0$.
\item[(b)] $Q$ degenerates so that the intrinsic distance between it opposite sides $L_2$ and $L_4$
tends to zero, while the intrinsic distance between its sides $L_1$ and $L_3$ does not tend to zero.
This happens when $\Q$ has an arc of order $1$ connecting $a_3$ with a point on $L_2$
mapped to $X_0$, but does not have an arc of order $1$ connecting $a_3$ to a point on $L_1$ mapped to $X_0$.
\item[(c)] $Q$ does not degenerate, but converges to a quadrilateral $Q'$ with the corner $a_3$ mapped to
the vertex $X_0$ of $\P$. This happens when $\Q$ does not have
an arc of order $1$ connecting $a_3$ with a point on one of the sides $L_1$ and $L_2$
mapped to $X_0$.
\end{itemize}

In the case (c), the quadrilateral $Q'$ must be one of the quadrilaterals
shown either in Fig.~\ref{3circles-em} or in Fig.~\ref{3circles-tetm},
since a quadrilateral $Q$ containing an irreducible quadrilateral $\tilde Q$
cannot have its integer corner $a_3$ mapped to a vertex of $\P$.
If $Q$ is one of the quadrilaterals in Fig.~\ref{3circles-em} then $Q'$ is equivalent to $Q$.
In this case, the combinatorial type of $Q$ does not change when $X$ passes trough $X_0$.
If $Q'$ is one of the quadrilaterals in Fig.~\ref{3circles-tetm} then it is not equivalent to $Q$,
and we say that $Q$ and $Q'$ are {\em adjacent}. In this case, $Q'$ has two adjacent quadrilaterals
with different combinatorial types.

\begin{figure}
\centering
\includegraphics[width=5.0in]{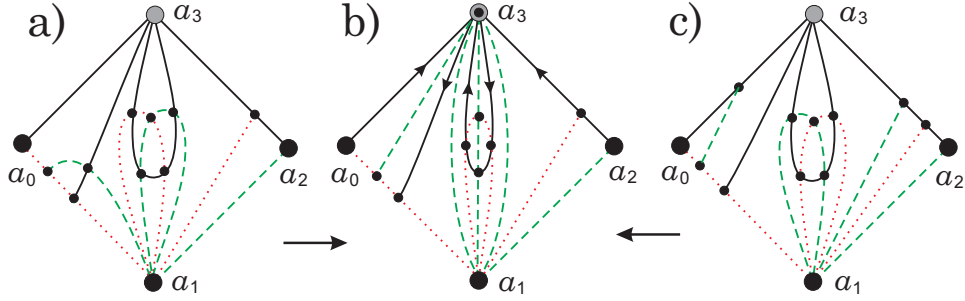}
\caption{An aa-chain of quadrilaterals of length 1.}\label{3circles-chain2}
\end{figure}

\begin{example}\label{ex:chain2}
The quadrilateral $Q$ in Fig.~\ref{3circles-chain2}a is $V^1_{1,1}$. As $X$ (the point on $C$
where $a_3$ maps) approaches $X_0\in C\cap C''$, the quadrilateral $Q$ converges to the
quadrilateral $Q'$ in Fig.~\ref{3circles-chain2}b which is the union of triangles $\bar M_1$ and $T_1$
with a digon $D_2$ inserted between them (see Fig.~\ref{3circles-tetm}d).
After $X$ passes $X_0$, the quadrilateral becomes $U^1_{1,1}$.
\end{example}

\begin{lemma}\label{adjacent} The following quadrilaterals are adjacent to the quadrilaterals in Fig.~\ref{3circles-tetm}:\newline
a) Let $Q'$ be the union of triangles $T_j$ and $E_n$ with a digon $D_{2\kappa}$ inserted between them
and digons $D_i$, $D_k$, $D_l$, $D_m$ attached as shown in Fig.~\ref{3circles-tetm}a.
Then $Q'$ is adjacent to the quadrilaterals $\bar X^\kappa_{2j,2n}$ and $\bar X^\kappa_{2j-1,2n+1}$
with digons $D_i$, $D_k$, $D_l$, $D_m$ attached as shown in Fig.~\ref{3circles-xruvz}b.\newline
b) Let $Q'$ be the union of triangles $T_j$ and $M_n$ with a digon $D_{2\kappa}$ inserted between them
and digons $D_i$, $D_k$, $D_l$, $D_m$ attached as shown in Fig.~\ref{3circles-tetm}b.
Then $Q'$ is adjacent to the quadrilaterals $\bar U^\kappa_{jn}$ and $\bar V^\kappa_{jn}$
with digons $D_i$, $D_k$, $D_l$, $D_m$ attached as shown in Fig.~\ref{3circles-xruvz}f,h.\newline
c) Let $Q'$ be the union of triangles $E_j$ and $T_n$ with a digon $D_{2\kappa}$ inserted between them
and digons $D_i$, $D_k$, $D_l$, $D_m$ attached as shown in Fig.~\ref{3circles-tetm}c.
Then $Q'$ is adjacent to the quadrilaterals $\bar X^\kappa_{2j+1,2n-1}$ and $\bar X^\kappa_{2j,2n}$
with digons $D_i$, $D_k$, $D_l$, $D_m$ attached as shown in Fig.~\ref{3circles-xruvz}b.\newline
d) Let $Q'$ be the union of triangles $\bar M_j$ and $T_n$ with a digon $D_{2\kappa}$ inserted between them
and digons $D_i$, $D_k$, $D_l$, $D_m$ attached as shown in Fig.~\ref{3circles-tetm}d.
Then $Q'$ is adjacent to the quadrilaterals $U^\kappa_{jn}$ and $V^\kappa_{jn}$
with digons $D_i$, $D_k$, $D_l$, $D_m$ attached as shown in Fig.~\ref{3circles-xruvz}e,g.
\end{lemma}

This can be shown by comparing the chains of the quadrilaterals in each of the cases a-d.

\begin{df}\label{df:chain}
{\rm For $k> 0$, a sequence $Q_0,Q'_1,Q_1,\dots,Q'_k,Q_k$ of
quadrilaterals with distinct combinatorial types,
where any two consecutive quadrilaterals are adjacent, and each of the terminal quadrilaterals
$Q_0$ and $Q_k$ has only one adjacent quadrilateral, is called a {\em chain} of the {\em length} $k$.
A quadrilateral $Q_0$ having no adjacent quadrilaterals
is called a chain of the length 0.

If both cases (a) and (b) are possible for degeneration of $Q_0$ and $Q_k$
then the chain is called an {\em ab-chain}.
If only the case (a) is possible, the chain is an {\em aa-chain}.
If only the case (b) is possible, the chain is a {\em bb-chain}.}
\end{df}

\begin{figure}
\centering
\includegraphics[width=5.0in]{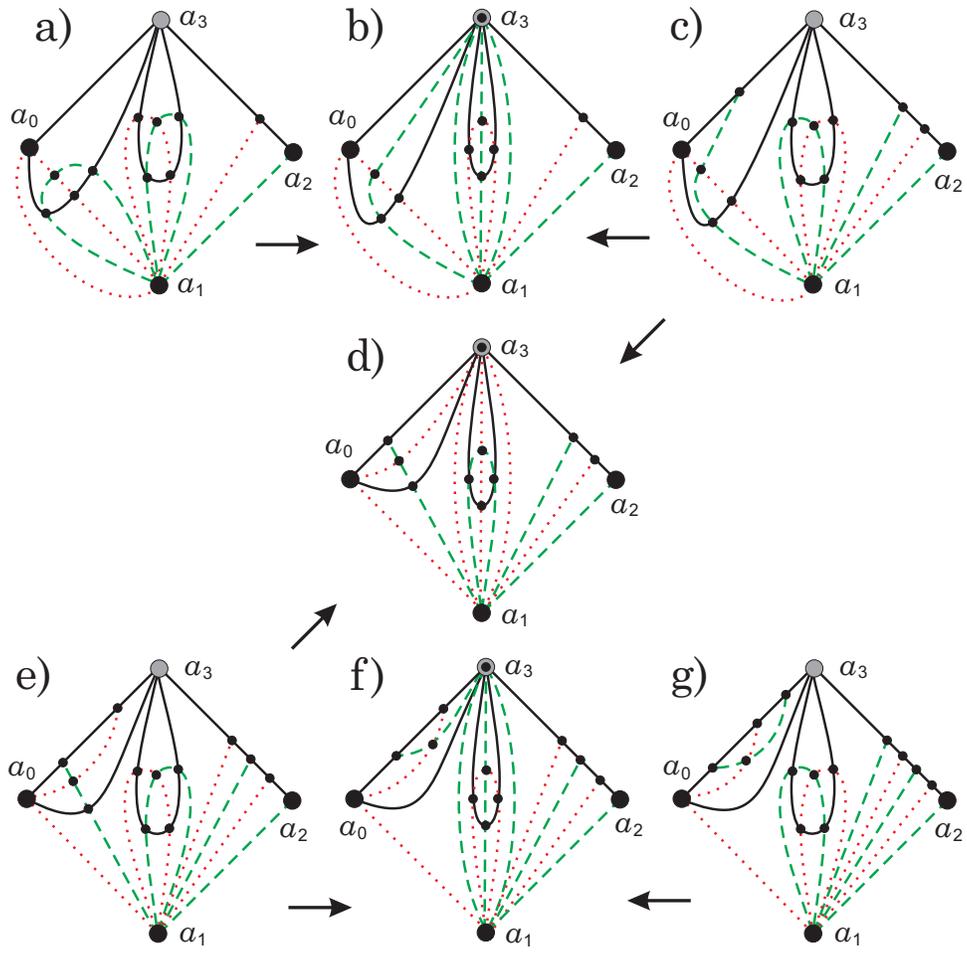}
\caption{An aa-chain of quadrilaterals of length 3.}\label{3circles-chain3}
\end{figure}

\begin{example}\label{ex:chain3}{\rm
A chain of quadrilaterals of length $3$ is shown in Fig.~\ref{3circles-chain3}.
The quadrilateral $Q_0$ in Fig.~\ref{3circles-chain3}a can be represented either as $V^1_{11}$ with a digon $D_1$ attached to
its side $L_1$ or as $\bar X^1_{3,1}$ with a digon $D_1$ attached to its side $L_4$.
The quadrilateral $Q_1$ in Fig.~\ref{3circles-chain3}c can be represented either as $U^1_{11}$ with
a digon $D_1$ attached to its side $L_1$ or as $\bar X^1_{2,2}$ with a digon $D_1$ attached to its side $L_4$.
The quadrilateral $Q_2$ in Fig.~\ref{3circles-chain3}e is $\bar X_{13}$
with a digon $D_1$ attached to its side $L_4$.
The quadrilateral $Q_3$ in Fig.~\ref{3circles-chain3}g is $\bar X_{04}$
with a digon $D_1$ attached to its side $L_4$.
The quadrilateral $Q'_1$ in Fig.~\ref{3circles-chain3}b is a union of triangles $\bar M_1$ and $T_0$ with
a digon $D_2$ inserted between them and a digon $D_1$ attached to its side $L_1$.
It can be also represented as a union of triangles $E_1$ and $T_0$ with
a digon $D_2$ inserted between them and a digon $D_1$ attached to its side $L_4$.
The quadrilateral $Q'_2$ in Fig.~\ref{3circles-chain3}d is a union of triangles $\bar T_0$ and $E_1$ with
a digon $D_2$ inserted between them and a digon $D_1$ attached to its side $L_4$.
The quadrilateral $Q'_3$ in Fig.~\ref{3circles-chain3}f is a union of triangles $E_0$ and $T_1$ with
a digon $D_2$ inserted between them and a digon $D_1$ attached to its side $L_4$.

If the point $X$ to which the corner $a_3$ of $Q_0$ maps approaches $C\cap C'$,
the intrinsic distance between the sides $L_1$ and $L_3$ (but not of $L_2$ and $L_4$) tends to zero (case a).
If $X$ approaches $C\cap C''$, the quadrilateral $Q_0$ converges to $Q'_1$ (case c).
If the point $X$ to which the corner $a_3$ of $Q_1$ maps approaches $C\cap C'$ (resp., $C\cap C''$),
 the quadrilateral $Q_1$ converges to $Q'_2$ (resp., $Q'_1$).
 If the point $X$ to which the corner $a_3$ of $Q_3$ maps approaches $C\cap C'$ (resp., $C\cap C''$),
 the quadrilateral $Q_2$ converges to $Q'_2$ (resp., $Q'_3$).
If the point $X$ to which the corner $a_3$ of $Q_3$ maps approaches $C\cap C''$,
 the quadrilateral $Q_3$ converges to $Q'_3$.
If $X$ approaches $C\cap C'$,
the intrinsic distance between the sides $L_1$ and $L_3$ (but not of $L_2$ and $L_4$) tends to zero (case a).
Thus the chain in Fig.~\ref{3circles-chain3} is an aa-chain.}
\end{example}

\begin{lemma}\label{aa-chains}
A chain of quadrilaterals with one integer angle
is an aa-chain (resp., a bb-chain) if and only if it contains
a quadrilateral $U^\kappa_{\mu\nu}$ (resp., $\bar U^\kappa_{\mu\nu}$) with
digons $D_i$, $D_k$, $D_l$, $D_m$ attached as in Fig.~\ref{3circles-xruvz}e,f,
with $\mu\ge 1$, $\nu\ge 1$, and $\min(i,l)=0$ (resp., $\min(k,m)=0$).
Each aa-chain (resp., bb-chain) contains only one quadrilateral
satisfying these conditions.
\end{lemma}

\noindent{\em Proof.}
Due to reflection symmetry, it is enough to consider only aa-chains
containing a quadrilateral $Q$ with the net $U^\kappa_{\mu\nu}$
and digons $D_i$, $D_k$, $D_l$, $D_m$ attached, as in Fig.~\ref{3circles-xruvz}e,
satisfying conditions of Lemma \ref{aa-chains}.

If $i=0$ then the net $\Q$ of $Q=Q_1$
contains an arc of length 1 connecting its corner $a_3$ with a point on its side $L_1$.
On the other hand, $Q$ is adjacent to a quadrilateral $Q'_1$ which is the union of
triangles $\bar M_\mu$ and $T_\nu$ with a digon $D_{2\kappa}$ inserted between them
and digons $D_k$, $D_l$, $D_m$ attached (since $i=0$, digon $D_i$ is empty)
as shown in Fig.~\ref{3circles-tetm}d (see Lemma \ref{adjacent}).
A quadrilateral $Q_0$ with the net $V^\kappa_{\mu\nu}$
and digons $D_k$, $D_l$, $D_m$ attached, as in Fig.~\ref{3circles-xruvz}g
is adjacent to $Q'_1$, and its net contains an arc of length 1 connecting its corner
$a_3$ with a point on its side $L_1$. Thus the chain of length 1 consisting of $Q_0$, $Q'_1$ and $Q_1=Q$
is an aa-chain. An example of such a chain is shown in Fig.~\ref{3circles-chain2}.
Note that if $l>0$ then $\nu=1$, and the quadrilateral $Q$ is equivalent to a quadrilateral with the
net $X^\kappa_{\mu,1}$ and digons $D_{k+1}$, $D_{l-1}$, $D_m$ attached to its sides $L_2$, $L_3$, $L_4$.

Suppose now that $i>0$ but $l=0$. This is only possible if $\mu=1$.
An example of such a net is shown in Fig.~\ref{3circles-chain3}c,
and it is shown in Example \ref{ex:chain3} that it belongs to an aa-chain.
The same arguments as in Example \ref{ex:chain3} show that any quadrilateral with the net
$U^\kappa_{1,\nu}$ and digons $D_i$, $D_k$, $D_m$ attached to its sides $L_1$, $L_2$, $L_4$,
where $i>0$, which is equivalent to a quadrilateral with the net $\bar X^\kappa_{2,2\nu}$ and
digons $D_{i-1}$, $D_k$, $D_{m+1}$ attached to its sides $L_1$, $L_2$, $L_4$,
belongs to an aa-chain of length 1, containing also a quadrilateral with the net $V^\kappa_{1,\nu}$
with the same digons as those attached to the sides of $U^\kappa_{1,\nu}$, and quadrilaterals
$\bar X^\kappa_{1,2\nu+1}$ $\bar X^\kappa_{0,2\nu+2}$ with the same digons as those attached
to the sides of $\bar X^\kappa_{2,2\nu}$.

It remains to show that any chain that either does not contain a quadrilateral $U^\kappa_{\mu\nu}$
with digons $D_i$, $D_k$, $D_l$, $D_m$ attached to its sides, or contains such a quadrilateral with
$i>0$ and $l>0$, is not an aa-chain. The corresponding statement for $\bar U^\kappa_{\mu\nu}$ would
follow by reflection symmetry.

We start with the chains containing $U^\kappa_{\mu\nu}$.
Conditions $i>0$ and $l>0$ are satisfied only when $\mu=\nu=1$.
Such a quadrilateral (for $\kappa=0$, $i=l=1$ and $k=m=0$) is shown in Fig.~\ref{3circle-chain1}c.
It is an easy exercise to show that it belongs to an ab-chain of length 2 containing also
$V^0_{1,1}$ with two digons $D_1$ attached to $L_1$ and $L_3$ (equivalent to $\bar X^0_{3,1}$
with two digons $D_1$ attached to $L_3$ and $L_4$, see Fig.~\ref{3circle-chain1}a) and $\bar V^0_{1,1}$ with two digons $D_1$
attached to $L_2$ and $L_4$ (equivalent to $\bar X^0_{1,3}$ with two digons $D_1$ attached to $L_3$ and $L_4$,
see Fig.~\ref{3circle-chain1}e).

The same arguments show that a quadrilateral $Q=Q_1$ with the net $U^\kappa_{1,1}$ and digons
$D_i$, $D_k$, $D_l$, $D_m$ attached to its sides (equivalent to a quadrilateral $X^\kappa_{2,2}$
with digons $D_i$, $D_{k+1}$, $D_{l-1}$, $D_m$ attached, to a quadrilateral $\bar X^\kappa_{2,2}$ with digons
$D_{i-1}$, $D_k$, $D_l$, $D_{m+1}$ attached, and to a quadrilateral  $\bar U^\kappa_{1,1}$
with digons $D_{i-1}$, $D_{k+1}$, $D_{l-1}$, $D_{m+1}$ attached) belongs to an ab-chain
of length 2 containing also a quadrilateral $Q_0$ with the net $V^\kappa_{1,1}$ and the same
digons as those attached to the sides of $U^\kappa_{1,1}$ (equivalent to
$\bar X^\kappa_{3,1}$ with the same digons as those attached to the sides of $\bar X^\kappa_{2,2}$),
and a quadrilateral $Q_2$ with the net $\bar V^\kappa_{1,1}$ and the same
digons as those attached to the sides of $\bar U^\kappa_{1,1}$ (equivalent to
$\bar X^\kappa_{1,3}$ with the same digons as those attached to the sides of $\bar X^\kappa_{2,2}$).

By reflection symmetry, all chains containing $\bar U^\kappa_{\mu\nu}$ are either bb-chains or ab-chains.
One can easily check that any chain containing a net $V^\kappa_{\mu\nu}$ (resp., $\bar V^\kappa_{\mu\nu}$) with digons
$D_i$, $D_k$, $D_l$, $D_m$ attached to its sides as in Fig.~\ref{3circles-xruvz}g,h contains also a net $U^\kappa_{\mu\nu}$
(resp., $\bar U^\kappa_{\mu\nu}$) with the same digons attached to its sides.

Now we are going to show that all chains that do not contain either $U^\kappa_{\mu\nu}$ or $\bar U^\kappa_{\mu\nu}$
are ab-chains.

Any net $X^\kappa_{\mu\nu}$ with digons $D_l$ and $D_m$ attached to its sides $L_3$ and $L_4$ has arcs of length 1
connecting $a_3$ with points on both $L_1$ and $L_2$, thus it is an ab-chain of length 0.
If, in addition, a non-empty digon $D_k$ is attached to its side $L_2$ then $\nu\le 1$.
When $\mu>0$, this net is equivalent to $U^\kappa_{\mu,1}$ with digons $D_{k-1}$, $D_{l+1}$, $D_m$ attached.
When $\mu=0$, it is equivalent to $\bar X^\kappa_{0,2}$ with digons $D_{k-1}$, $D_{l+1}$, $D_m$ attached,
and belongs to an ab-chain of length 2 containing also the nets $\bar X^\kappa_{1,1}$ and $\bar X^\kappa_{2,0}$
with the same digons as those attached to the sides of $\bar X^\kappa_{0,2}$.
Similarly, the net $X^\kappa_{0,1}$ with digons $D_i$, $D_k$, $D_l$, $D_m$ attached, $k>0$,
is equivalent to $\bar X^\kappa_{0,2}$ with digons $D_i$, $D_{k-1}$, $D_{l+1}$, $D_m$ attached,
and belongs to an ab-chain of length 2 containing also the nets $\bar X^\kappa_{1,1}$ and $\bar X^\kappa_{2,0}$
with the same digons as those attached to the sides of $\bar X^\kappa_{0,2}$.

All nets $\bar X^\kappa_{\mu\nu}$ with digons $D_i$, $D_k$, $D_l$, $D_m$ attached to their sides
as in Fig.~\ref{3circles-xruvz}b have been already considered, except those with $l=m=0$.
For given $\kappa$, $\mu+\nu$, $i$ and $k$, such nets belong to a single ab-chain of length $\mu+\nu$.

Each net $R_{\mu,\nu}$ (resp., $\bar R_{\mu,\nu}$) with digons $D_i$, $D_k$, $D_l$, $D_m$ attached
to its sides as in Fig.~\ref{3circles-xruvz}c,d has arcs of length 1
connecting $a_3$ with points on both $L_1$ and $L_2$, thus it is an ab-chain of length 0.
The same is true for each net $Z^\kappa_{\mu\nu}$ with digons $D_i$, $D_k$, $D_l$, $D_m$ attached
to its sides as in Fig.~\ref{3circles-xruvz}i.

Each net in Fig.~\ref{3circles-em} is an ab-chain of length 0.
Although its net may not have arcs of length 1
connecting $a_3$ with points on both $L_1$ and $L_2$,
the net does not change combinatorial equivalence class when the image $X$
of its corner $a_3$ passes through an intersection point of $C$ with either $C'$ or $C''$,
unless the arc of $\nabla$ connecting $a_3$ with either $a_0$ or $a_1$ contracts to a point.

This completes the proof of Lemma \ref{aa-chains}.

\begin{lemma}\label{u-chains}
For the given angles $\alpha_0,\dots,\alpha_3$ at the corners $a_0,\dots,a_3$ of a quadrilateral
$Q$ with one integer corner $a_3$,
the number of nets $U^\kappa_{\mu\nu}$ with digons $D_i$, $D_k$, $D_l$, $D_m$ attached
to its sides as in Fig.~\ref{3circles-xruvz}e, so that $\min(i,l)=0$, equals
\begin{equation}\label{u-count}
\left[\min\left(\frac{\alpha_3}2,\frac{1+[\alpha_1]}2,\frac{\delta}2\right)\right]
\end{equation}
where $[\alpha_0]=i+m$, $[\alpha_1]=i+k+2\kappa+\nu$, $[\alpha_2]=k+l$, $\alpha_3=l+m+2\kappa+\mu+1$
are integer parts of the angles of $Q$, and
$\delta=\max(0,1+[\alpha_1]+\alpha_3-[\alpha_0]-[\alpha_2])/2$.
\end{lemma}

\noindent{\em Proof.}
One can establish one-to-one correspondence between the nets in Lemma \ref{u-chains}
and the corresponding nets in Theorem 17.2 of \cite{EGT2} for a quadrilateral with two
opposite integer corners, except the angle at the corner $a_2$ of the quadrilateral
in Theorem 17.2 of \cite{EGT2} would be $1+[\alpha_1]$.

\begin{thm}\label{thm:count}
For given angles $\alpha_0,\dots,\alpha_3$, where $\alpha_3$ is integer
and $\alpha_0,\dots,\alpha_2$ are not integer, and for a fixed modulus,
there are at least
\begin{equation}\label{main-count}
a_3-2\left[\min\left(\frac{\alpha_3}2,\frac{1+[\alpha_1]}2,\frac{\delta}2\right)\right]
\end{equation}
quadrilaterals, where $\delta=\max(0,1+[\alpha_1]+\alpha_3-[\alpha_0]-[\alpha_2])/2$.
\end{thm}

\noindent{\em Proof.}
It follows from Proposition \ref{allreal} that the total number of chains of
quadrilaterals with the angles as in Theorem \ref{thm:count} is $a_3$.
It follows from Lemmas \ref{aa-chains} and \ref{u-chains}
that the number of aa-chains (and the number of bb-chains which is equal to
the number of aa-chains by reflection symmetry) is given by (\ref{u-count}).
Thus the number of ab-chains is given by (\ref{main-count}).
For each ab-chain there is a quadrilateral with such chain with any modulus,
and these quadrilaterals are distinct for distinct chains.

\bigskip
{\em A. E. and A. G.: Department of Mathematics, Purdue University,

West Lafayette, IN 47907-2067 USA
\vspace{.1in}

V. T.: Department of Mathematical Sciences, IUPUI,

Indianapolis, IN 46202-3216 USA;

St. Petersburg branch of Steklov Mathematical Institute,

Fontanka 27, St.Petersburg, 191023, Russia.}


\bigskip
\begin{thebibliography}{11}
\bibitem{Ahlfors} L.~Ahlfors, Conformal invariants. Topics in geometric
function theory, Reprint of the 1973 original, AMS Chelsea Publishing,
Providence, RI, 2010.
\bibitem{E} A.~Eremenko, Metrics of positive curvature with conic
singularities on the sphere, Proc.~Amer.~Math.~Soc., 132 (2004) 3349--3355.
\bibitem{EGT1} A.~Eremenko, A.~Gabrielov and V.~Tarasov, Metrics with conic
singularities and spherical polygons, arXiv:1405.1738.
\bibitem{EGT2} A.~Eremenko A.~Gabrielov and V.~Tarasov, Metrics with four conic
singularities and spherical quadrilaterals, arXiv:1409.1529.
\bibitem{FFKRUY} S.~Fujimori, Y.~Kawakami, M.~Kokubu, W.~Rossman,
M.~Umehara and K.~Yamada, CMC-1 trinoids in hyperbolic 3-space and metrics
of constant curvature one with conical singularities on the 2-sphere,
Proc.~Japan Acad., 87 (2011) 144--149.
\bibitem{Krein} F. Gantmakher and M. Krein, Oscillation matrices and kernels
and small vibrations of mechanical systems, AMS, Chelsea, Providence, RI,
2011.
\bibitem{Klein1} F.~Klein, Mathematical seminar at G\"ottingen,
winter semester 1905/6 under the direction of Professors Klein,
Hilbert and Minkowski, talks by F.~Klein, notes by O.~Toeplitz,
\newline
www.claymath.org/sites/default/files/klein1math.sem{\_}.{\_}ws1905-06.pdf
\bibitem{Klein-book} F.~Klein, Forlesungen \"uber die hypergeometrische
Funktion, reprint of the 1933 original, Springer Verlag, Berlin-NY,
1981.
\bibitem{Pont} L. Pontrjagin,
Hermitian operators in spaces with indefinite metric, Bull. Acad. Sci. URSS,
S\'er. Math. (Izvestia Akad. Nauk SSSR) 8, (1944) 243--280.
\bibitem{Schil} F.~Schilling, Ueber die Theorie der symmetrischen
s-Funktion mit einem einfachen Nebenpunkte,  Math.~Ann., 51, (1899) 481-522.
\bibitem{Troy1} M.~Troyanov, Metrics of constant curvaturer
on a sphere with two conical singularities,
Lect.~Notes Math., 1410, Springer, Berlin, 1989, 296--306.
\end{thebibliography}
\end{document}